\documentclass{amsart}
\usepackage[english,francais]{babel}
\usepackage{amsmath}

\usepackage{amsfonts}

\usepackage{hyperref}

\usepackage[ansinew]{inputenc}
\font\bb=msbm7 at 11pt
\font\mathbb=msbm7 at 8pt

\def \F {\hbox{\bb F}}
\def \N {\hbox{\bb N}}
\def \R {\hbox{\bb R}}

\def \A {\hbox{\bb A}}
\def \Z {\hbox{\bb Z}}
\def \G {\hbox{\bb G}}
\def \T {\hbox{\bb G}}
\def \p {\mathfrak p}
\def \z {\hbox{\mathbb Z}}
\def \n {\hbox{\mathbb N}}
\def \g {\hbox{\mathbb G}}
\def \Q {\hbox{\bb Q}}
\def \AA {\mathcal{A}}
\def \L {\mathcal{L}}
\def \M {\mathcal{M}}
\def \FF {\mathcal{F}}
\def \U {\mathcal{U}}
\def \H {\mathcal{H}}
\def \P {\hbox{\mathbb P}}
\def \x {\mbox{\bf x}}
\def \y {\mbox{\bf y}}
\def \r {\mbox{\bf r}}
\def \ii {\mbox{\bf i}}
\def \k {\mbox{\bf k}}
\def \u {\mbox{\bf u}}

\def \e {\mbox{\bf e}}
\def \f {\mbox{\bf f}}
\def \r {\mbox{\bf r}}
\def \s {\mbox{\bf s}}

\def \MM {\mbox{\bf M}}
\def \GL {\mbox{\bf GL}}
\newcommand{\dd}{\boldsymbol{\delta}}
\newcommand{\ee}{\boldsymbol{\varepsilon}}
\newcommand{\Ker}{\ensuremath{\mbox{\rm{Ker}}}}
\newcommand{\Ima}{\ensuremath{\mbox{\rm{Im}}}}
\newcommand{\Vect}{\ensuremath{\mbox{\rm{Vect}}}}
\newcommand{\NP}{\ensuremath{\mbox{\rm{NP}}}}
\newcommand{\GNP}{\ensuremath{\mbox{\rm{GNP}}}}

\newcommand{\HP}{\ensuremath{\mbox{\rm{HP}}}}
\newcommand{\HS}{\ensuremath{\mbox{\rm{HS}}}}

\newcommand{\Tr}{\ensuremath{\mbox{\rm{Tr}}}}
\newcommand{\pr}{\ensuremath{\mbox{\rm{pr}}}}
\newcommand{\ppcm}{\ensuremath{\mbox{\rm{ppcm}}}}
\def \ep {\varepsilon}

\newtheorem{theo}{Th\'eor\`eme}[section]
\newtheorem{lem}{Lemme}[section]
\newtheorem{defi}{D\'efinition}[section]
\newtheorem{coro}{Corollaire}[section]
\newtheorem{propo}{Proposition}[section]

\newtheorem{rem}{Remarque}[section]

\author{R\'egis Blache}
\address{Equipe AOC,
IUFM de la Guadeloupe}
\email{blache@iufm.univ-ag.fr}

\title[Polygones de Newton de certaines sommes de caractères.]{Polygones de Newton de certaines sommes de caractères et séries de Poincaré.}
\begin{document}

\begin{abstract}
On se propose dans cet article de donner quelques r\'esultats sur le comportement asymptotique des polygones de Newton des fonctions $L$ associ\'ees \`a des sommes exponentielles, provenant de certains polyn\^omes de Laurent en $n$ variables. A cet effet, on étudie et on utilise la somme directe de polyèdres convexes. Cette opération permet de d\'eterminer ais\'ement la limite des polygones de Newton g\'en\'eriques associés à la somme directe $\Delta=\Delta_1\oplus \Delta_2$ quand on conna\^{\i}t la limite des polygones de Newton génériques associés à chacun des polyèdres $\Delta_i$. Ce sont à notre connaissance les premiers résultats sur le comportement asymptotique des polygones de Newton pour les polynômes à plusieurs variables. 
\end{abstract}

%\begin{altabstract}
%In this paper, we shall precise the asymptotic behaviour of Newton polygons of $L$ functions associated to character sums, coming from some $n$ variable Laurent polynomials. In order to do this, we use the free sum on convex polytopes. This operation allows the determination of the limit of generic Newton polygons for the sum $\Delta=\Delta_1\oplus \Delta_2$ when we know the limit of generic Newton polygons for each factor. To our knowledge, these are the first results concerning the asymptotic behaviour of Newton polygons for multivariable polynomials.
%\end{altabstract}

\subjclass{11M38,14F30,52B20}
\keywords{Sommes de caract\`eres, fonctions $L$, polygones et polyèdres de Newton}
%\altkeywords{Character sums, $L$-functions, Newton polygons and polytopes}
\maketitle

\tableofcontents

\bigskip

\section{Introduction}

\bigskip

Dans toute la suite, on note $k$ un corps fini poss\'edant $q=p^a$ \'el\'ements, et $k_r$ son extension de degr\'e $r$ dans une cl\^oture alg\'ebrique $\overline{k}$ fix\'ee une fois pour toutes. On note $\x=(x_1,\dots,x_n)$, et on choisit $f(\x)=\sum_{\ii \in\z^n}a_{\ii}\x^{\ii}\in k[\x,\x^{-1}]$ un polyn\^ome de Laurent en $n$ variables \`a coefficients dans $k$. Si $\psi$ est un caract\`ere additif (non trivial) de $k$, on note $\psi_r:=\psi\circ\Tr_{k_r/k}$ le caract\`ere induit par $\psi$ sur $k_r$ ; d'autre part on note $\chi$ un caract\`ere multiplicatif de $(k^\times)^n$, et $\chi_r:=\chi\circ \mbox{\rm N}_{k_r/k}$ son extension \`a $(k_r^\times)^n$. On forme les sommes exponentielles associ\'ees \`a $f$ et $\chi$ sur chacune des extensions de $k$
$$S_r(f,\chi)=\sum_{\x\in \g_m^n(k_r)} \psi(f(\x))\chi(\x),$$
puis \`a partir de cette famille de sommes on construit la fonction $L$
$$L(f,\chi;T)=\exp\left(\sum_{r\geq 1} S_r(f,\chi)\frac{T^r}{r}\right).$$
Quand $\chi$ est le caract\`ere trivial, on notera simplement $L(f;T)$ cette fonction. On sait depuis les travaux de Dwork et de Grothendieck ({\it cf.} \cite{Dw}, \cite{Gr}) que c'est une fraction rationnelle. 

\medskip

Commen\c{c}ons par le cas o\`u $\chi$ est trivial, qui est le plus classique dans la litt\'erature. Le premier résultat dans cette direction est dû à Deligne \cite[Théorème 8.4]{De}. Il affirme que pour un polynôme de degré $d$ premier à $p$ dont la forme de plus haut degré définit une hypersurface non singulière de l'espace projectif $\P^{n-1}$, la fonction $L'(f,T)$ (où cette fois les sommes sont définies sur $\A^n$ et non plus sur $\G_m^n$) est de degré $(d-1)^n$.

De façon générale, on peut associer au polyn\^ome $f$ son {\it polyèdre de Newton \`a l'infini}, qui est le polyèdre convexe $\Delta$ d\'efini dans l'espace affine $\R^n$ comme l'enveloppe convexe de l'origine et du support de $f$ (c'est \`a dire de l'ensemble des $\ii$ de $\Z^n$ tels que $a_{\ii}$ est non nul). Alors sous des conditions de non d\'eg\'en\'erescence de certaines formes de plus haut degr\'e de $f$, Adolphson et Sperber \cite{AS1} ont montr\'e que $L(f;T)^{(-1)^{n-1}}$ est en fait un polyn\^ome de degr\'e $n!V(\Delta)$, o\`u $V$ d\'esigne le volume usuel sur $\R^n$. 

Notons $\alpha_1,\dots,\alpha_{n!V(\Delta)}$ les racines réciproques de ce polynôme. Ce sont des entiers algébriques possédant les propriétés suivantes : si $\alpha_i$ est l'un d'entre eux, son module complexe est $|\alpha_i|=q^{\frac{w_i}{2}}$, pour un entier $0\leq w_i \leq n$. Tous les conjugués de $\alpha_i$ sont de même module complexe. La distribution des $w_i$ est connue \cite[Theorem 1.8]{DL}. D'autre part, pour tout premier $\ell\neq p$, $\alpha_i$ est une unité $\ell$-adique. Finalement, on a $|\alpha_i|_p=q^{-s_i}$, pour $s_i$ un rationnel compris entre $0$ et $n$. C'est à ces valuations que nous allons nous intéresser. 

Comme d'habitude, on identifie les valuations $q$-adiques des racines r\'eciproques d'un polyn\^ome avec les pentes des segments de son polygone de Newton $q$-adique. Pour simplifier les notations, nous noterons $\NP_q(f)$ ({\it resp.} $\NP_q(f,\chi)$) le polygone de Newton $q$-adique de $L(f;T)$ ({\it resp.} de $L(f,\chi;T)$) dans la suite. Si $\Pi$ est un polygone convexe de longueur $l$, c'est à dire le graphe d'une fonction convexe, continue sur l'intervalle $[0,l]$ et affine par morceaux sur chacun des intervalles $[i-1,i]$, on notera $\Pi=(\pi_i)_{1\leq i \leq l}$ quand sa pente sur $[i-1,i]$ est $\pi_i$. Si $\Pi_1$ et $\Pi_2$ sont deux polygones convexes de longueur $l$, on écrira $\Pi_1 \preceq \Pi_2$ quand $\Pi_1$ est au dessus de $\Pi_2$, et que leurs extrémités coïncident

Adolphson et Sperber ont prouv\'e \cite[Theorem 3.10]{AS1} l'existence d'une borne inf\'erieure pour les polygones de Newton des fonctions $L(f;T)^{(-1)^{n-1}}$ quand $f$ d\'ecrit l'ensemble des polyn\^omes non d\'eg\'en\'er\'es de polyèdre $\Delta$. Cette borne est souvent appel\'ee {\it polygone de Hodge}, et not\'ee $\HP(\Delta)$. Il s'agit d'un invariant ne d\'ependant que de $\Delta$, que nous allons maintenant d\'ecrire. Remarquons qu'il ne s'agit en g\'en\'eral ni du polygone de Hodge d'un cristal comme dans \cite{Ka}, ni d'un polygone de Hodge g\'eom\'etrique, ces deux familles de polygones ayant toutes leurs pentes enti\`eres.

Notons $C(\Delta):=\R_+\Delta$ le c\^one de $\Delta$ dans $\R^n$, $M_\Delta:=C(\Delta)\cap\Z^n$ le mono\"{\i}de associ\'e \`a ce c\^one, et $\AA_{\Delta}$ l'alg\`ebre $k[\x^{M_\Delta}]$. On peut d\'efinir une application de $C(\Delta)$ dans $\R_+$, le {\it poids associ\'e \`a $\Delta$} par
$$w_{\Delta}(\u)=\min\{\rho\in \R_+,~\u\in\rho\Delta\}.$$
Les sommets du polyèdre $\Delta$ \'etant dans $\Z^n$, l'image de $M_\Delta$ par $w_\Delta$ est contenue dans $\Q_+$, et plus pr\'ecis\'ement dans $\frac{1}{D}\N$ pour un certain entier $D>0$, minimal, qu'on appellera dans la suite le {\it d\'enominateur de $\Delta$}. Le poids $w_\Delta$ fait de l'alg\`ebre $\AA_\Delta$ une alg\`ebre gradu\'ee
$$\AA_\Delta=\oplus_{i\geq 0} \AA_{\Delta,\frac{i}{D}},~\AA_{\Delta,\frac{i}{D}}=\Vect\{\x^{\u},~w_\Delta(\u)=\frac{i}{D}\}$$
\`a laquelle on associe sa s\'erie de Poincar\'e
$$P_{\AA_\Delta}(t):=\sum_{i\geq 0} \dim\AA_{\Delta,\frac{i}{D}}t^i.$$
Kouchnirenko \cite[Lemme 2.9]{Kou} a montr\'e que quand $f$ est non d\'eg\'en\'er\'e, cette s\'erie est en fait une fraction rationnelle. Plus pr\'ecis\'ement, $P_\Delta(t):=(1-t^D)^n P_{\AA_\Delta}(t)$ est un polyn\^ome de degr\'e inférieur ou égal à $nD$. Si on note $P_\Delta(t):=\sum \ell_i t^{s_i}$, le polygone $\HP(\Delta)$ est alors le polygone convexe commen\c{c}ant \`a l'origine, et form\'e de la juxtaposition des segments de longueur horizontale $\ell_i$ et de pente $\frac{s_i}{D}$. Par commodité, on l'appellera aussi dans la suite le {\it polygone issu de la série de Poincaré $P_{\AA_\Delta}$}. Le r\'esultat d'Adolphson et Sperber se r\'e\'ecrit donc ainsi : pour tout polyn\^ome $f$ de $k[\x,\x^{-1}]$ de polyèdre $\Delta$ et non d\'eg\'en\'er\'e, on a $\NP_q(f)\preceq \HP(\Delta)$.

Il est maintenant naturel  de se demander comment varient les polygones $\NP_q(f)$ quand $f$ varie parmi les polyn\^omes de polyèdre fix\'e, non d\'eg\'en\'er\'es. Malheureusement ces variations sont tr\`es difficiles \`a contr\^oler ; des calculs explicites dans le cas de polyn\^omes de petit degr\'e en une variable montrent qu'il est illusoire d'esp\'erer donner une r\'eponse compl\`ete \`a cette question. Pour contourner cette difficult\'e, on pr\'ef\`ere parler de {\it polygone de Newton g\'en\'erique}. Le th\'eor\`eme de sp\'ecialisation de Grothendieck \cite{Ka} assure que la borne inf\'erieure des $\NP_q(f)$ existe, et qu'elle est atteinte pour tous les points dans un ouvert dense de l'espace des polyn\^omes de polyèdre fix\'e, non d\'eg\'en\'er\'es ; c'est cette borne inf\'erieure qu'on appelle le polygone de Newton g\'en\'erique. Ce polygone ne d\'epend pas de $q$, mais seulement de $p$, et on le note $\GNP(\Delta,p)$. Dans le cas de la dimension $1$, on sait calculer explicitement ce polygone \cite{BF}, ainsi que le polyn\^ome de Hasse, qui d\'etermine l'hypersurface de l'espace des polyn\^omes hors de laquelle on a $\NP_q(f)=\GNP(\Delta,p)$.

On se pose donc la question du comportement du polygone de Newton g\'en\'erique. Pour des raisons li\'ees \`a la ramification, il est ais\'e de voir qu'une condition n\'ecessaire pour qu'il co\"{\i}ncide avec le polygone de Hodge est d'avoir $p\equiv 1\mod D$. Adolphson et Sperber ont conjectur\'e que cette condition est suffisante. La conjecture est av\'er\'ee en dimension $n\leq 3$, mais elle est fausse en dimension sup\'erieure, o\`u il faut remplacer $D$ par un multiple $D^*$ en g\'en\'eral strict, comme l'a d\'emontr\'e Wan \cite{Wan1}, \cite{Wan2}. Donc on a $\liminf_{p\rightarrow \infty} \GNP(\Delta,p)=\HP(\Delta)$. Wan a conjectur\'e \cite[Conjecture 1.11]{Wan2} que la limite existe sous certaines hypothèses, c'est \`a dire que $\lim_{p\rightarrow \infty} \GNP(\Delta,p)=\HP(\Delta)$. Ce résultat est connu pour les polynômes de Laurent en une variable \cite{BF}, \cite{Liu}. On va d\'emontrer cette conjecture dans deux cas 
\begin{itemize}
	\item[\it i/] quand $\Delta$ est l'enveloppe convexe de points de la forme $\{d_i\e_i,-d_i'\e_i\}_{1\leq i\leq n}$, avec $(\e_1,\dots,\e_n)$ une base du $\Z$-module $\Z^n$ ; ce sera le th\'eor\`eme \ref{lim}. 
	\item[\it ii/] quand $\Delta$ est l'enveloppe convexe de points de la forme $\{d_i\f_i,-d_i'\f_i\}_{1\leq i\leq n}$, où $\f_1,\dots,\f_n$ engendrent un sous module $M$ de $\Z^n$ tel que $2\Z^n\subset M$ ; ce sera le th\'eor\`eme \ref{lim2}.
\end{itemize}

Ce sont à notre connaissance les premiers résultats sur le comportement asymptotique de somme de caractères associées à des polynômes de plusieurs variables. Notons que le cas {\it i/} couvre en particulier les polynômes étudiés par Deligne.

\medskip

Une autre question, plus difficile, est la suivante : choisissons un polyn\^ome de Laurent $\widetilde{f}$ à coefficients dans $\overline{\Q}$, et soit $\Q_{\widetilde{f}}$ l'extension de $\Q$ engendrée par les coefficients de $\widetilde{f}$. Pour chaque premier $p$ de $\Q$, on choisit $\p$ un premier au dessus de $p$ dans le corps $\Q_f$, de corps résiduel $\F_q$. On se demande comment varient les polygones de Newton $\NP_q(\widetilde{f}\mod \p)$ des r\'eductions modulo $\p$ de $\widetilde{f}$ quand $p$ tend vers l'infini. Considérons l'espace des polyn\^omes à coefficients dans $\overline{\Q}$, de polyèdre de Newton $\Delta$, et de monomes prescrits de façon à ce que le sous-monoïde de $M_{\Delta}$ engendré par les exposants des monômes prescrits contienne tous les éléments de $M_{\Delta}$ à l'exception d'un nombre fini. Alors Wan conjecture \cite[Conjecture 1.12]{Wan2} qu'il existe un ouvert dense défini sur $\Q$ de l'espace de ces polyn\^omes tel que pour tout $\widetilde{f}$ de cet ouvert on ait 
$$\lim_{p\rightarrow \infty}\NP_q(\widetilde{f}\mod \p)=\HP(\Delta).$$ 
Ce r\'esultat est connu pour l'espace de tous les polynômes de degré $d$ en une variable \cite[Theorem 1.3]{Zhu1}, ainsi que pour l'espace des polynômes de Laurent de degrés $d$ et $d'$ en une variable \cite{LZ}. On va le montrer pour certains espaces de polyn\^omes \`a plusieurs variables dont le polyèdre de Newton est de l'une des formes décrites dans les cas {\it i/} et {\it ii/} ci-dessus ; voir les théorèmes \ref{dense} et \ref{dense2}.

\medskip

Dans le cas où $\chi$ n'est pas trivial, la situation est un peu plus compliquée. Adolphson et Sperber ont montré, toujours sous des conditions de non-dégénérescence, que la fonction $L$ a le même degré que dans le cas additif ({\it cf.} \cite{AS3}, \cite{AS2}). Ils ont aussi donné une borne inférieure pour les polygones de Newton de telles sommes. Si on peut toujours décrire cette borne inférieure à l'aide de séries de Poincaré, elle dépend maintenant du résidu de $p$ modulo l'ordre du caractère $\chi$. En particulier on ne peut plus espérer obtenir une limite quand $p$ tend vers $\infty$, à moins de supposer que le caractère est d'ordre deux (c'est ce résultat qui nous permet de prouver le cas {\it ii/} ci-dessus). En revanche, quand $p$ tend vers $\infty$ le long d'une classe modulo l'ordre du caractère, on retrouve l'existence d'une limite, généralisant ainsi les résultats en dimension $1$ de \cite{BFZ}. Ces résultats font l'objet des théorèmes \ref{limS} et \ref{denseS}.

\medskip

Cet article est organisé de la façon suivante : dans le premier chapitre, nous définissons la somme directe de polyèdres convexes, et calculons le polygone de Hodge associé. Dans le second nous utilisons des résultats de cohomologie $\ell$-adique (principalement la formule de Künneth) pour exprimer les polygones de Newton de (fonctions $L$ de) sommes associées à certains polynômes en plusieurs variables à l'aide des polygones de Newton associés à des polynômes plus simples. Dans le chapitre 3 sont démontrés les conjectures de Wan dans le cas {\it i/} : on y rappelle la situation, connue, des polynômes de Laurent en une variable, puis on en déduit les théorèmes \ref{lim} et \ref{dense}. Ensuite, dans le quatrième chapitre, on donne quelques applications au cas des sommes mixtes, tordues par un caractère multiplicatif. Finalement, on utilise les deux résultats et la formule de Poisson pour déduire les conjectures de Wan pour le cas {\it ii/} dans le dernier chapitre.

\section{Sommes directes de polyèdres.}

Dans toute cette section, on fixe deux polyèdres convexes $\Delta_1$ de $\R^{n_1}$ et $\Delta_2$ de $\R^{n_2}$. On va rappeler la définition de leur {\it somme directe} $\Delta_1\oplus \Delta_2$, et donner certaines de ses propriétés ; on exprimera ensuite la s\'erie de Poincar\'e de l'alg\`ebre gradu\'ee $\AA_{\Delta_1\oplus \Delta_2}$ \`a l'aide de celles des alg\`ebres gradu\'ees $\AA_{\Delta_1}$ et $\AA_{\Delta_2}$, de façon à calculer le polygone de Hodge de $\Delta_1\oplus \Delta_2$ en fonction de ceux de $\Delta_1$ et $\Delta_2$.

\medskip

Commençons par définir la somme directe de polyèdres ({\it cf.} \cite[16.1.3]{hrz}) ; dans toute la suite, on va supposer que les polyèdres qui interviennent sont de dimension maximale, c'est à dire que l'espace affine qu'ils engendrent est l'espace ambiant.

\begin{defi}
\label{somdir}
Soient deux polyèdres convexes $\Delta_1$ et $\Delta_2$, respectivement dans $\R^{n_1}$ et $\R^{n_2}$. Leur {\it somme directe} est le polyèdre convexe de $\R^{n_1+n_2}$ qui est l'enveloppe convexe de $\Delta_1\times \{0\}\cup \{0\}\times\Delta_2$. On le note $\Delta_1\oplus\Delta_2$.
\end{defi}

\begin{rem}
Il ne faut pas confondre l'op\'eration qu'on vient de d\'efinir avec la somme usuelle de polyèdres (ou somme de Minkowski). Par exemple, si $\Delta_1=[0,d_1]\subset \R$ et $\Delta_2=[0,d_2]\subset \R$, alors $\Delta_1\oplus\Delta_2$ est le triangle de sommets $(0,0),(d_1,0)$ et $(0,d_2)$, alors que la somme usuelle de ces polyèdres donne le rectangle de sommets $(0,0),(d_1,0),(0,d_2)$ et $(d_1,d_2)$.
\end{rem}

On va déterminer les faces ne contenant pas l'origine de la somme directe de deux polyèdres convexes contenant chacun l'origine. Ces résultats sont sans doute bien connus mais comme de nombreuses démonstrations de cet article en dépendent, et pour n'avoir pas trouvé de référence convenable, nous les rappelons ici. Par définition, les faces d'un polyèdre sont ses intersections avec ses hyperplans d'appui, ainsi que le polyèdre lui même et l'ensemble vide $\emptyset$ ; ces deux dernières sont parfois appelées impropres. 

\begin{propo}
\label{faces}
Soient $\Delta_1$ et $\Delta_2$ deux polyèdres convexes contenant l'origine, et $\Delta:=\Delta_1\oplus\Delta_2$ leur somme directe. Les faces de $\Delta$ ne contenant pas l'origine sont les polyèdres $\sigma:=\sigma_1\oplus \sigma_2$, où $\sigma_i$ décrit les faces de $\Delta_i$ ne contenant pas l'origine pour $i=1,2$.

Si de plus $\sigma_i$ est une face de dimension $d_i$, alors $\sigma$ est une face de $\Delta$ de dimension $d_1+d_2+1$.

\end{propo}

\begin{proof}
Commençons par démontrer que si $\sigma_1$ ({\it resp.} $\sigma_2$) est une face (éventuellement vide) de $\Delta_1$ ({\it resp.} $\Delta_2$) ne contenant pas l'origine, alors $\sigma_1\oplus \sigma_2$ est une face de $\Delta$ ne contenant pas l'origine. Soit $H_1$ d'équation $\sum_{i=1}^{n_1} a_ix_i=1$ ({\it resp.} $H_2$ d'équation $\sum_{i=1}^{n_2} b_iy_i=1$) un hyperplan d'appui de $\Delta_1$ ({\it resp.} $\Delta_2$) pour sa face $\sigma_1$ ({\it resp.} $\sigma_2$). Pour une face vide, on prendra tous les coefficients nuls. Alors puisque $\Delta_1$ contient l'origine, il est contenu dans le demi espace $H_1^-=\{(x_1,\dots,x_{n_1}),~\sum_{i=1}^{n_1} a_ix_i\leq 1\}$ de $\R^{n_1}$, et il en est de même pour $\Delta_2$ dans $\R^{n_2}$. Considérons l'hyperplan $H$ de $\R^{n_1+n_2}$ d'équation $\sum_{i=1}^{n_1} a_ix_i+\sum_{i=1}^{n_2} b_iy_i=1$. Alors par définition de $\Delta$ comme enveloppe convexe, ce polyèdre est contenu dans le demi espace $H^-$ de $\R^{n_1+n_2}$, et l'intersection $H\cap \Delta$ contient par construction $\sigma_1\times\{0\}$ et $\{0\}\times\sigma_2$. En particulier $H$ est un hyperplan d'appui de $\Delta$, et la face $\sigma:=H\cap\Delta$ qu'il détermine contient l'enveloppe convexe de $\sigma_1\times\{0\}\cup \{0\}\times\sigma_2$, c'est à dire la somme directe $\sigma_1\oplus \sigma_2$. Pour montrer l'autre inclusion, choisissons $\z(z_1,\dots,z_{n_2})$ un point de $\sigma$. Alors $\z$ est dans $\Delta$, barycentre de $(\x,0)$ et $(0,\y)$ deux points respectivement dans $\Delta_1\times\{0\}$ et $\{0\}\times \Delta_2$ ; on peut donc trouver un réel $\lambda\in[0,1]$ tel que $z_i=\lambda x_i$ pour $1\leq i\leq n_1$ et $z_{n_1+i}=(1-\lambda)y_i$ pour $1\leq i\leq n_2$. Mais puisque $\z$ est dans $H$, on a $\sum a_iz_i+\sum b_iz_{n_1+i}=1$ ; d'autre part $\sum_{i=1}^{n_1} a_ix_i\leq 1$ et $\sum_{i=1}^{n_2} b_iy_i\leq 1$. Ces deux dernières inégalités doivent être des égalités, c'est à dire que $\x\in \sigma_1$, $\y\in \sigma_2$, et $x\in \sigma_1\oplus \sigma_2$. On a donc démontré que $\sigma=\sigma_1\oplus \sigma_2$ est une face de $\Delta$ ne contenant pas l'origine.

\medskip

Inversement soit $\sigma$ une face de $\Delta$ ne contenant pas l'origine, et $H$ un hyperplan d'appui pour cette face. Puisque $\sigma$ est un polyèdre convexe, il est l'enveloppe convexe de ses points extrémaux, qui sont les points extrémaux de $\Delta$ contenus dans $\sigma$. Mais par construction les points extrémaux de $\Delta$ sont de la forme $(\x_1,0)$ ou $(0,\x_2)$, où $\x_i$ décrit les points extrémaux de $\Delta_i$. Notons $S_1$ l'ensemble des points extrémaux de $\sigma$ du premier type, $S_2$ l'ensemble des points extrémaux de $\sigma$ du second. Alors $\sigma$ est l'enveloppe convexe de $S_1\cup S_2$, c'est à dire la somme directe des convexes $\sigma_1$ et $\sigma_2$ avec $\sigma_i$ l'enveloppe convexe de $S_i$, et il suffit maintenant de montrer que $\sigma_i$ est une face de $\Delta_i$. Si $S_i=\emptyset$, il n'y a rien à montrer ; sinon par construction $H_i=H\cap\R^{n_i}$ est un hyperplan d'appui de $\Delta_i$ dans $\R^{n_i}$, et $\sigma_i=\Delta\cap H_i$ est bien une face de $\Delta_i$.

\medskip

Prouvons finalement l'assertion sur la dimension : la face $\sigma_i$ est l'enveloppe convexe d'un ensemble $S_i$ de points ne contenant pas l'origine parmi lesquels on peut choisir un sous ensemble maximal $S_i'$ de $d_i+1$ points affinement indépendants. Les ensembles de points de $\R^{n_1+n_2}$, $T_1=S_1'\times\{0\}$ et $T_2=\{0\}\times S_2'$ sont disjoints, et affinement indépendants par construction. L'ensemble $T=T_1\cup T_2$ est donc un sous ensemble maximal de $d_1+d_2$ points affinement indépendants de $S= S_1\times\{0\}\cup \{0\}\times S_2$ ; mais comme $\sigma$ est par construction l'enveloppe convexe de $S$, l'assertion en résulte.
\end{proof}
 
Comme $\Delta_i$ contient l'origine, ses faces ne contenant pas l'origine sont nécessairement de dimension inférieure ou égale à $n_i-1$ ; on déduit donc le
 
\begin{coro}
\label{codim1}
 Les faces de codimension $1$ de $\Delta$ ne contenant pas l'origine sont les $\sigma=\sigma_1\oplus \sigma_2$, où $\sigma_i$ décrit les faces de codimension $1$ de $\Delta_i$ ne contenant pas l'origine.
 
 De plus si $H_1$ d'équation $\sum a_ix_i=1$ ({\it resp.} $H_2$ d'équation $\sum b_iy_i=1$) est l'hyperplan d'appui de $\Delta_1$ pour $\sigma_1$ ({\it resp.} de $\Delta_2$ pour $\sigma_2$), alors l'hyperplan d'appui de $\Delta$ pour $\sigma$ a pour \'equation $\sum a_ix_i+\sum b_iy_i=1$.
 \end{coro}

On va maintenant exprimer les diff\'erents objets associ\'es au polyèdre $\Delta_1\oplus\Delta_2$ dans l'introduction \`a l'aide de ceux associ\'es \`a chacun de ses facteurs. On aura besoin d'une nouvelle description du poids : soit $\u(u_1,\dots,u_{n_1+n_2})$ un point de $C(\Delta)$. La demi droite $\R^+\u$ rencontre la frontière du polyèdre $\Delta$ en un point d'une face $\sigma$ de codimension $1$ ne contenant pas l'origine. Si $H$ est l'hyperplan d'appui pour cette face, d'équation $a_1x_1+\dots+a_{n_1+n_2}x_{n_1+n_2}=1$, alors le poids de $\u$ est $$w_{\Delta}(\u)=a_1u_1+\dots+a_{n_1+n_2}u_{n_1+n_2}.$$

\begin{lem}
Soient $\Delta_1$ et $\Delta_2$ deux polyèdres convexes contenant l'origine, et $\Delta:=\Delta_1\oplus\Delta_2$ leur somme directe. On note $\sigma_i$ une face de $\Delta_i$. Alors
\begin{itemize}
	\item[i/] le c\^one $C(\sigma_1\oplus\sigma_2)$ dans $\R^{n_1+n_2}$ est égal au produit $C(\sigma_1)\times C(\sigma_2)$ ;
	\item[ii/] le mono\"{\i}de $M_{\Delta}$ est égal au mono\"{\i}de $M_{\Delta_1}\times M_{\Delta_2}$ de $\Z^{n_1+n_2}$ ;
	\item[iii/] le poids $w_{\Delta}$ est l'application $w_{\Delta_1}+w_{\Delta_2}$ de $C(\Delta)$ dans $\R_+$ qui \`a $\u=(\u_1,\u_2)$ associe $w_{\Delta_1}(\u_1)+w_{\Delta_2}(\u_2)$
  \item[iv/] le d\'enominateur $D$ de $\Delta$ est le plus petit commun multiple des d\'enominateurs $D_1$ et $D_2$ de $\Delta_1$ et $\Delta_2$.
\end{itemize}
\end{lem}

\begin{proof}
Montrons la première assertion : le point $\u(\u_1,\u_2)$ est dans le cone $C(\sigma_1\oplus\sigma_2)$ si et seulement si on peut trouver un réel $\rho$ tel que $(\rho\u_1,\rho\u_2)$ soit dans $\sigma_1\oplus\sigma_2$. Mais par définition, ce dernier polyèdre convexe est formé des points de la forme $(\lambda\x_1,(1-\lambda)\x_2)$ quand $\lambda$ d\'ecrit $[0,1]$ et les $\x_i$ d\'ecrivent $\sigma_i$. Chaque $\u_i$ est donc dans $C(\sigma_i)$, et la réciproque est évidente.

Le second point est maintenant une conséquence facile de la définition du monoïde associé à un polyèdre convexe, et de l'assertion {\it i/} appliquée \`a $\sigma_i=\Delta_i$.

Pour montrer {\it iii/}, on va utiliser le corollaire \ref{codim1}, ainsi que l'expression du poids donnée ci-dessus. Notons $\sigma$ une face de codimension $1$ de $\Delta$ ne contenant pas l'origine en un point de laquelle la demi droite $\R_+ \u$ rencontre la frontière de $\Delta$.  D'après le corollaire \ref{codim1}, on a $\sigma=\sigma_1\oplus\sigma_2$, pour $\sigma_i$ une face de codimension $1$ de $\Delta_i$ ne contenant pas l'origine. En particulier si $\u=(\u_1,\u_2)$, le {\it i/} nous assure que $\u_i\in C(\sigma_i)$. Si on note $H_1 : a_1x_1+\dots+a_{n_1}x_{n_1}=1$ ({\it resp.} $H_2 : a_{n_1+1}x_{n_1+1}+\dots+a_{n_1+n_2}x_{n_1+n_2}=1$) l'hyperplan d'appui de $\Delta_1$ pour $\sigma_1$ ({\it resp.} de $\Delta_2$ pour $\sigma_2$) et si $\u_1(u_1,\dots,u_{n_1})$ ({\it resp.} $\u_2(u_{n_1+1},\dots,u_{n_1+n_2})$), on doit donc avoir $w_{\Delta_1}(\u_1)=a_1u_1+\dots+a_{n_1}u_{n_1}$ ({\it resp. } $w_{\Delta_2}(\u_2)=a_{n_1+1}u_{n_1+1}+\dots+a_{n_1+n_2}u_{n_1+n_2}$). Toujours d'après le corollaire \ref{codim1}, l'équation de l'hyperplan d'appui de $\Delta$ pour $\sigma$ est $a_1x_1+\dots+a_{n_1+n_2}x_{n_1+n_2}=1$, et $w_{\Delta}(\u)=a_{1}u_{1}+\dots+a_{n_1+n_2}u_{n_1+n_2}$ ; c'est le résultat annoncé. 

La derni\`ere assertion est une cons\'equence directe de {\it iii/} et de la d\'efinition du d\'enominateur d'un polyèdre convexe. 
\end{proof}

On d\'eduit de l'assertion {\it ii/} que l'alg\`ebre $\AA_\Delta$ est isomorphe au produit tensoriel (sur $k$) des alg\`ebres $\AA_{\Delta_1}$ et $\AA_{\Delta_2}$. Venons en \`a la graduation ; si $\x_i^{\u_i}$ est dans $\AA_{\Delta_i,\frac{k_i}{D_i}}$, alors d'apr\`es l'assertion {\it iii/}, le mon\^ome $\x^{\u}=\x_1^{\u_1}\x_2^{\u_2}$ est dans $\AA_{\Delta,\frac{k}{D}}$, avec 
$$\frac{k_1}{D_1}+\frac{k_2}{D_2}=\frac{k}{D}.$$ 
On obtient la d\'ecomposition suivante pour chaque pi\`ece de la graduation de $\AA_\Delta$ :
$$\AA_{\Delta,\frac{k}{D}}=\bigoplus_{\frac{k_1}{D_1}+\frac{k_2}{D_2}=\frac{k}{D}}\AA_{\Delta_1,\frac{k_1}{D_1}}\otimes \AA_{\Delta,\frac{k_2}{D_2}},$$
puis la factorisation de la s\'erie de Poincar\'e de $\AA_\Delta$ \`a l'aide de celles de $\AA_{\Delta_1}$ et $\AA_{\Delta_2}$
$$P_{\AA_\Delta}(t)=P_{\AA_{\Delta_1}}(t^{\frac{D}{D_1}})P_{\AA_{\Delta_2}}(t^{\frac{D}{D_2}}),$$
et finalement la factorisation $P_{\Delta}(t)=P_{\Delta_1}(t^{\frac{D}{D_1}})P_{\Delta_2}(t^{\frac{D}{D_2}})$. 

\medskip

Nous terminons ce paragraphe en d\'emontrant, \`a l'aide de la formule pr\'ec\'edente, que le polygone de Hodge de la somme directe $\Delta=\Delta_1\oplus\Delta_2$ s'exprime \`a l'aide des polygones de Hodge de chacun des facteurs. Pour cela nous avons besoin d'introduire une nouvelle op\'eration sur les polygones convexes. Rappelons qu'on a choisi de noter un polygone convexe d'origine $O$ par $(s_i)_{1\leq i\leq a}$ quand il est form\'e par la juxtaposition des segments de longueur horizontale $1$ et de pente $s_i$.

\begin{defi}
Soient deux polygones convexes $\Pi_1$ et $\Pi_2$. Alors si 
$$\Pi_1=(s_i)_{1\leq i\leq a},~\Pi_2=(s_i')_{1\leq i\leq b},$$ 
on d\'efinit le produit de $\Pi_1$ et $\Pi_2$, et on note $\Pi_1\times \Pi_2$ le polygone convexe d'origine $O$ d\'efini par
$$\Pi=(s_i+s_j')_{1\leq i\leq a,~1\leq j\leq b}.$$
\end{defi}

Remarquons que la longueur horizontale de $\Pi$ est le produit des longueurs horizontales de $\Pi_1$ et $\Pi_2$, mais aussi que la longueur horizontale du segment de pente $s$ dans $\Pi$ est 
$$\ell=\sum_{s_i+s_j'=s}\ell_i\ell_j',$$
où $\ell_i$ ({\it resp.} $\ell_j'$) est la longueur horizontale du segment de pente $s_i$ ({\it resp.} $s_j'$) de $\Pi_1$ ({\it resp.} $\Pi_2$).

On en déduit la d\'ecomposition suivante pour le polygone $\HP(\Delta)$.

\begin{propo}
\label{prodH}
Soient $\Delta_1$ et $\Delta_2$ deux polyèdres convexes, et $\Delta$ leur somme directe. Alors le polygone de Hodge de $\Delta$ est le produit des polygones de Hodge de ses facteurs
$$\HP(\Delta)=\HP(\Delta_1)\times \HP(\Delta_2).$$
\end{propo}

\begin{proof}
Soit $\ell^{(1)}_{k_1}$ ({\it resp.} $\ell^{(2)}_{k_2}$) la longueur horizontale du segment de $\HP(\Delta_1)$ ({\it resp.} $HP(\Delta_2)$) de pente $\frac{k_1}{D_1}$ ({\it resp.} $\frac{k_2}{D_2}$). Le segment de pente $\frac{k}{D}$ du produit $\HP(\Delta_1)\times \HP(\Delta_2)$ a pour longueur horizontale $\ell=\sum\ell^{(1)}_{k_1}\ell^{(2)}_{k_2}$ o\`u la somme porte sur les $k_1$, $k_2$ tels que $\frac{k_1}{D_1}+\frac{k_2}{D_2}=\frac{k}{D}$. Mais la construction de $\HP(\Delta)$ \`a partir de $P_\Delta$, jointe \`a la factorisation de $P_\Delta$, donne le m\^eme r\'esultat.
\end{proof}

\section{Sommes exponentielles.}

Ici on consid\`ere deux polyn\^omes de Laurent sur $k$, $f_1$ et $f_2$ respectivement en $n_1$ et $n_2$ variables, d'ind\'etermin\'ees $\x_1$ et $\x_2$, et on appelle $f=(f_1,f_2)$ le polyn\^ome de Laurent en les $n:=n_1+n_2$ variables $\x=(x_1,\dots,x_{n_1+n_2})$ tel que $f(\x)=f_1(x_1,\dots,x_{n_1})+f_2(x_{n_1+1},\dots,x_{n_1+n_2})$. Il résulte immédiatement de la définition \ref{somdir} que si $\Delta_1$ et $\Delta_2$ sont les polyèdres respectifs de $f_1$ et $f_2$ dans $\R^{n_1}$ et $\R^{n_2}$, alors le polyèdre de $f$ est $\Delta$, la somme directe de $\Delta_1$ et $\Delta_2$.

On va dans ce chapitre r\'eexprimer les espaces vectoriels de cohomologie $\ell$-adique associ\'es aux sommes exponentielles provenant de $f$ \`a l'aide de la formule de K\"unneth, et de ceux associ\'es \`a $f_1$ et $f_2$ ; ensuite on en déduira des bornes pour les polygones de Newton génériques associés à $\Delta$ en fonction de ceux associés à $\Delta_1$ et $\Delta_2$. 

Commençons par montrer que la condition de non dégénérescence se transmet de $f_1$ et $f_2$ à $f$.

\begin{lem}
\label{nondeg}
Soient $f_1$ et $f_2$ deux polynômes non dégénérés respectivement pour $\Delta_1$ et $\Delta_2$. Alors le polynôme $f=(f_1,f_2)$ est non dégénéré pour $\Delta$.
\end{lem}

\begin{proof}
Rappelons que si $f$ est un polynôme de polyèdre de Newton $\Delta$, et $\sigma$ une face du polyèdre $\Delta$, le polynôme $f_\sigma$ est la somme des monômes de $f$ dont le support est dans $\sigma$. Alors $f$ est non dégénéré quand pour toute face de $\Delta$ ne contenant pas l'origine, les polynômes $\frac{\partial f_\sigma}{\partial x_i}$, $1\leq i\leq n_1+n_2$ n'ont pas de zéro commun dans $(\bar{k}^\times)^{n_1+n_2}$. Mais d'après le lemme \ref{faces}, toute face de $\Delta$ ne contenant pas l'origine est de la forme $\sigma_1\oplus\sigma_2$, pour $\sigma_1$ une face de $\Delta_1$ ({\it resp.} $\sigma_2$ une face de $\Delta_2$). Il est facile de vérifier que $f_\sigma(\x)=f_{1,\sigma_1}(\x_1)+f_{2,\sigma_2}(\x_2)$, et qu'on a
$$\frac{\partial f_\sigma}{\partial x_i}=\left\{\begin{array}{lll}
\frac{\partial f_{1,\sigma_1}}{\partial x_i} & \mbox{\rm pour } & 1\leq i\leq n_1 \\
\frac{\partial f_{2,\sigma_2}}{\partial x_i} & \mbox{\rm pour } & n_1+1\leq i\leq n_1+n_2 \\
\end{array} \right.$$
Donc $\x=(\x_1,\x_2)$ est un zéro commun des $\frac{\partial f_\sigma}{\partial x_i}$, $1\leq i\leq n_1+n_2$ quand $\x_1$ est un zéro commun des $\frac{\partial f_{1,\sigma_1}}{\partial x_i}$, $1\leq i\leq n_1$ et $\x_2$ est un zéro commun des $\frac{\partial f_{2,\sigma_2}}{\partial x_i}$, $n_1+1\leq i\leq n_1+n_2$. C'est à dire que la non dégénérescence de $f_1$ et $f_2$ assure celle de $f$.
\end{proof}

Soit $\psi$ un caract\`ere additif non trivial sur $k$, et $\L_\psi$ le $\overline{\Q}_\ell$-faisceau sur $\A^1_k$ associ\'e \`a $\psi$ et au recouvrement d'Artin-Schreier $y^q-y=x$. De même on note $\chi$ un caractère de $k^\times$ et $\L_{\chi}$ le $\overline{\Q}_\ell$-faisceau sur $\G_{m,k}$ associ\'e. 

Si $X$ est un sch\'ema de type fini sur $k$, $f$ une fonction r\'eguli\`ere sur $X$ (c'est \`a dire un morphisme $f:X\rightarrow \A^1$), et $g$ une fonction régulière qui ne s'annule pas sur $X$, on peut construire comme dans l'introduction la fonction $L(X,f,g;T)$, et la formule des traces de Grothendieck nous permet de la r\'einterpr\'eter \`a l'aide des polyn\^omes caract\'eristiques de l'action du Frobenius sur les groupes de cohomologie du faisceau $f^*\L_\psi\otimes g^*\L_\chi$
$$L(X,f,g;T)=\prod_i \det\left(I-TF|H^i_c(X\otimes \overline{k},f^*\L_\psi\otimes g^*\L_\chi)\right)^{(-1)^{i-1}}.$$
 
Revenons \`a la situation qui nous int\'eresse. On a ici les trois fonctions $f_i:\T_m^{n_i}\rightarrow \A^1$, $1\leq i\leq 2$, et $f=(f_1,f_2):\T_m^{n}\rightarrow \A^1$. On fixe un caractère $\chi_1$ ({\it resp.} $\chi_2$) de $(k^\times)^{n_1}$ ({\it resp.} $(k^\times)^{n_2}$), et on note $\chi$ le caractère $(\chi_1,\chi_2)$ de $(k^\times)^{n}$. D'apr\`es la d\'efinition de $f$, on a, en notant $\pr_i$ les projections canoniques de $\T_m^n=\T_m^{n_1}\times \T_m^{n_2}$ sur chacun de ses facteurs, que $f^*\L_\psi\otimes \L_{\chi}=\bigotimes_{i=1}^2 \pr_i^*\left(f_i^*\L_\psi\otimes \L_{\chi_i}\right)$ est le produit tensoriel externe des deux faisceaux $f_i^*\L_\psi\otimes \L_{\chi_i}$. Alors la formule de K\"unneth nous assure qu'on a un isomophisme 
$$H^\bullet_c(\T_m^n,f^*\L_\psi\otimes \L_{\chi})=H^\bullet_c(\T_m^{n_1},f_1^*\L_\psi\otimes \L_{\chi_1})\otimes H^\bullet_c(\T_m^{n_2},f_2^*\L_\psi\otimes \L_{\chi_2}).$$
On a vu que $f$ est non d\'eg\'en\'er\'e quand $f_1$ et $f_2$ le sont ; dans ce cas l'isomorphisme pr\'ec\'edent se r\'e\'ecrit simplement
$$H^n_c(\T_m^n,f^*\L_\psi\otimes \L_{\chi})=H^{n_1}_c(\T_m^{n_1},f_1^*\L_\psi\otimes \L_{\chi_1})\otimes H^{n_2}_c(\T_m^{n_2},f_2^*\L_\psi\otimes \L_{\chi_2})$$
d'apr\`es des r\'esultats de Denef et Loeser \cite[Theorem 1.3]{DL}. En d'autres termes, la fonction $L(f,\chi;T)^{(-1)^{n-1}}$ est le polyn\^ome dont les racines r\'eciproques sont les produits des couples de racines r\'eciproques des polyn\^omes $L(f_1,\chi_1;T)^{(-1)^{n_1-1}}$ et $L(f_2,\chi_2;T)^{(-1)^{n_2-1}}$. 

Rappelons que pour un polyn\^ome de Laurent $f$ on note $\NP_q(f,\chi)$ le polygone de Newton du polyn\^ome $L(f,\chi;T)^{(-1)^{n-1}}$. On d\'eduit en particulier des r\'esultats pr\'ec\'edents une factorisation de $\NP_q(f,\chi)$ qui nous sera utile un peu plus loin. 

\begin{lem}
\label{prodnp}
On a l'\'egalit\'e de polygones de Newton 
$$\NP_q(f,\chi)=\NP_q(f_1,\chi_1)\times NP(f_2,\chi_2).$$
\end{lem}

Pour un polyèdre $\Delta$ de dimension $n$, et $\chi$ un caractère multiplicatif comme ci-dessus, on a défini le polygone de Newton g\'en\'erique $\GNP(\Delta,\chi,p)$ comme la borne inf\'erieure des polygones de Newton $\NP_q(f,\chi)$ quand $f$ parcourt l'ensemble des polyn\^omes de polyèdre $\Delta$, non d\'eg\'en\'er\'es. Quand $\Delta$ est une somme directe, on peut déduire des résultats précédents une borne pour ce polygone à l'aide des polygones de Newton génériques des facteurs de $\Delta$.

\begin{coro}
\label{prodgnp}
Soient $\Delta_1$ et $\Delta_2$ deux polyèdres convexes, et $\Delta$ leur somme directe. Alors on a  
$$\GNP(\Delta_1,\chi_1,p)\times \GNP(\Delta_2,\chi_2,p)\preceq\GNP(\Delta,\chi,p).$$
\end{coro}

\begin{proof}
Le th\'eor\`eme de sp\'ecialisation de Grothendieck (voir par exemple \cite{Ka}) nous assure que pour chaque $i$, il existe un ouvert dense $\U_{\Delta_i,\chi_i,p}$ de $\M_{\Delta_i}$, l'espace des coefficients des polyn\^omes de polyèdre $\Delta_i$, non d\'eg\'en\'er\'es, tel que pour tout $f_i$ de $\U_{\Delta_i,\chi_i,p}$, on ait $\NP_q(f_i,\chi_i)=\GNP(\Delta_i,\chi_i,p)$. D'après le lemme \ref{nondeg}, si  $f_1$ et $f_2$ sont tous deux non d\'eg\'en\'er\'es, alors $f=(f_1,f_2)$ l'est aussi pour $\Delta$, donc on a l'inclusion $\M_{\Delta_1}\times\M_{\Delta_2}\subset \M_\Delta$. On d\'eduit du lemme \ref{prodnp} que la borne inf\'erieure des polygones de Newton $\NP_q(f,\chi)$ quand $f$ parcourt $\M_{\Delta_1}\times\M_{\Delta_2}$ est le polygone $\GNP(\Delta_1,\chi_1,p)\times \GNP(\Delta_2,\chi_2,p)$ (atteinte pour tous les polynômes à coefficients dans $\U_{\Delta_1,\chi_1,p}\times\U_{\Delta_2,\chi_2,p}$). Mais d'apr\`es la d\'efinition de $\GNP(\Delta,\chi,p)$ comme borne inf\'erieure des $\NP_q(f,\chi)$ quand $f$ d\'ecrit $\M_\Delta$, on  obtient le r\'esultat.
\end{proof}

\section{Comportement asymptotique, cas additif.}

On se place dans la situation suivante : on fixe un entier $n$, et on note $(\e_1,\dots,\e_n)$ une base du $\Z$-module $\Z^n$. D'autre part on choisit des entiers naturels $d_1,d_1',\dots,d_n,d_n'$ tels que pour chaque $i$ on ait $(d_i,d_i')\neq (0,0)$ (sinon la situation qu'on va d\'ecrire peut se ramener en dimension inf\'erieure). On note $\Delta$ le polyèdre convexe de $\R^n$ qui est l'enveloppe convexe de l'ensemble des points $\{d_i\e_i,-d_i'\e_i\}_{1\leq i\leq n}$ et de l'origine si n\'ecessaire.

On se propose de d\'emontrer le r\'esultat suivant :

\begin{theo}
\label{lim}
Quand $p$ tend vers l'infini, le polygone de Newton g\'en\'erique de $\Delta$ associ\'e au premier $p$, $\GNP(\Delta,p)$, tend vers le polygone de Hodge $\HP(\Delta)$.
\end{theo}

La d\'emonstration, qu'on effectuera un peu plus loin, est une cons\'equence des r\'esultats pr\'ec\'edents et d'un cas d\'ej\`a connu de ce th\'eor\`eme, en dimension $1$. Commen\c{c}ons par rappeler ce qu'on sait en dimension $1$ ; le lecteur int\'eress\'e par les d\'etails pourra se r\'ef\'erer aux articles \cite{BF}, ,\cite{Liu}, \cite{Zhu1}.

Soit $f$ un polyn\^ome de Laurent en la variable $x$, $f(x)=\sum_{i=-d'}^d a_ix^i$, $a_{-d'}a_d\neq 0$. Il est clair que le polyèdre convexe associ\'e \`a $f$ est le segment de $\R$ d'extr\'emit\'es $-d'$ et $d$ ; le poids est donné, pour tout $n\in \Z$, par $w(n)=\max(\frac{n}{d},-\frac{n}{d'})$. On d\'eduit de la définition que le polygone de Hodge $\HP([-d',d])$ est le polygone d'extr\'emit\'es l'origine et le point de coordonn\'ees $(d+d',\frac{d+d'}{2})$, et poss\'edant un segment de longueur $1$ pour chacune des pentes suivantes
$$0,1,\frac{1}{d},\dots,\frac{d-1}{d},\frac{1}{d'},\dots,\frac{d'-1}{d'} ~\left( 0,\frac{1}{d},\dots,\frac{d-1}{d} ~\mbox{\rm si}~d'=0\right).$$ 
Nous noterons d\'esormais $s_1,\dots,s_{d+d'}$ ces pentes, rang\'ees par ordre croissant. On en déduit en particulier une autre description du polygone $\HP([-d',d])$ : c'est le polygone issu de l'origine et passant par les points de coordonnées $(i,s_1+\dots+s_i)$ pour $1\leq i\leq d+d'$.

Comme dans \cite{Ro}, on sait associer \`a $f$ un op\'erateur diff\'erentiel sur un espace de s\'eries surconvergentes, ainsi qu'un op\'erateur de Frobenius, qui commutent, de fa\c{c}on \`a r\'einterpr\'eter la fonction $L(f;T)$ comme le polyn\^ome caract\'eristique de l'op\'erateur de Frobenius agissant sur le premier espace de cohomologie de de Rham. On peut alors estimer, pour $p$ assez grand, les parties principales des coefficients de cette matrice, et donner des congruences pour ses mineurs principaux, qui sont les coefficients de la fonction $L$. 

Notons $\pi$ l'unique racine du polyn\^ome $X^{p-1}+p$ dans une cl\^oture alg\'ebrique fixée du corps $\Q_p$ des nombres $p$-adiques telle que $\psi(1)\equiv 1+\pi [\pi^2]$. On a alors $\Q_p(\pi)=\Q_p(\zeta_p)$, et on pose $K=\Q_p(\zeta_p,\zeta_{q-1})$. On note, pour chaque \'el\'ement $a$ de $k$, $\widetilde{a}$ son rel\`evement de Teichm\"uller (si $a=0$, alors $\widetilde{a}=0$, sinon la r\'eduction de $\widetilde{a}$ modulo l'id\'eal maximal est $a$, et on a $\widetilde{a}\in \boldsymbol{\mu}_{q-1}$). Soit $\widetilde{f}$ le polyn\^ome de $K[x,x^{-1}]$ obtenu \`a partir de $f$ en relevant ses coefficients comme ci-dessus. Si on pose $L(f;T)=1+\sum_{i=1}^{d+d'} M_iT^i$, on obtient pour tout $1\leq i\leq d+d'$ la congruence (dans $K$)
$$M_i\equiv u_i\P_{d,d',i}^{\rho,\rho'}(\widetilde{a}_{-d'},\dots, \widetilde{a}_d)\pi^{aY_i}\mod\pi^{aY_i+1},$$
o\`u $u_i$ est une unit\'e de l'anneau de valuation de $K$, et les polyn\^omes $\P_{d,d',i}^{\rho,\rho'}$ peuvent \^etre choisis \`a coefficients dans $\Q$, ne d\'ependant que des degr\'es $d$ et $d'$, et des restes $\rho$ et $\rho'$ respectifs des divisions euclidiennes de $p$ par $d$ et $d'$. 

Si $0\leq i_1\leq d$ et $0\leq i_2< d'$ sont les deux entiers déterminés (de façon unique) par la condition 
$$\{ s_1,\dots,s_i\}=\{0\}\cup \{\frac{j}{d},~1\leq j\leq i_1\}\cup\{\frac{j}{d'},~1\leq j\leq i_2\},$$ 
alors on sait exprimer $Y_i$ \`a l'aide du groupe sym\'etrique sur $i$ \'el\'ements, agissant sur l'ensemble $\{-i_2,\dots,0,\dots,i_1\}$ :
$$Y_i=\min_{\sigma\in S_i} \sum_{j=-i_2}^{i_1} \lceil w(pj-\sigma(j))\rceil,$$
o\`u $\lceil x\rceil$ d\'esigne le plus petit entier sup\'erieur ou \'egal \`a $x$. 

On obtient donc une description assez pr\'ecise de la situation en dimension $1$ pour $p$ fixé et choisi assez grand : d'une part le polygone de Newton g\'en\'erique (par rapport à la valuation $v_q$), $\GNP([-d',d],p)$, qui a pour sommets l'origine et les $(i,\frac{Y_i}{p-1})_{1\leq i\leq d+d'}$, et d'autre part le polyn\^ome de Hasse 
$$\H_{[-d',d]}^{\rho,\rho'}(\widetilde{a}_{-d'},\dots, \widetilde{a}_d):=\prod_{i=1}^{d+d'} \P_{d,d',i}^{\rho,\rho'}(\widetilde{a}_{-d'},\dots, \widetilde{a}_d),$$
qui d\'efinit une hypersurface de l'espace des polyn\^omes de Laurent de degr\'es $d,d'$ à coefficients dans $\F_q$. Tous les polyn\^omes dont les coefficients ne sont pas dans cette hypersurface satisfont $\NP_q(f)=\GNP(d,d',p)$. 

Quand $p$ varie, il est aisé de vérifier que pour tout $1\leq i\leq d+d'$, on a 
$$\lim_{p\rightarrow \infty} \frac{Y_i}{p-1}=s_1+\dots+s_i,$$ 
c'est à dire que le polygone de Newton générique converge vers le polygone de Hodge. Finalement, définissons le polynôme (à coefficients dans $\Q$)
$$\H_{[-d',d]}(X_{-d'},\dots,X_d)=\prod_{(\rho,\rho') \in (\z/d\z)^\times\times (\z/d'\z)^\times}\H_{[-d',d]}^{\rho,\rho'}(X_{-d'},\dots, X_d).$$
De plus si $f\in\overline{\Q}[x,x^{-1}]$ est un polynôme de Laurent de la forme $f(x)=\sum_{i=-d'}^d A_ix^i$, et si pour tout premier $p$ on choisit $\p$ un premier au dessus de $p$ dans le corps $\Q_f$ engendré par les coefficients de $f$, on a $\lim_{p\rightarrow\infty} \NP_q(f\mod \p)=\HP([-d',d])$ dès que
$$\H_{[-d',d]}(A_{-d'},\dots,A_d)\neq 0,$$
c'est à dire qu'il existe un ouvert dense $\U$, défini sur $\Q$, de l'espace des polynômes de polyèdre $[-d',d]$ en une variable à coefficients dans $\overline{\Q}$, tel que pour tout $f$ dans $U$, la limite $\lim_{p\rightarrow\infty} \NP_q(f\mod \p)$ existe et coïncide avec le polygone de Hodge.

\medskip

Nous sommes maintenant à même de démontrer le théorème \ref{lim}.

\begin{proof} Commen\c{c}ons par r\'eduire le probl\`eme au cas o\`u $\e_i$ est le $i$-\`eme vecteur de la base canonique de $\R^n$. On peut d\'efinir une action à gauche de $\MM_n(\Z)$ sur l'espace des polyn\^omes de Laurent en $n$ variables, qui au polyn\^ome $f(\x)=\sum a_{\ii} \x^{\ii}$ et \`a la matrice $M$ associe le polynôme $^Mf(\x)=f(^M\x)=\sum a_{\ii} (^M\x)^{\ii}=\sum a_{\ii} \x^{M\ii}$, o\`u $^M\x$ est le $n$-uplet de variables dont la $i$\`eme est $\prod_{j=1}^n x_j^{m_{ji}}$, et $M\ii$ désigne la multiplication usuelle de la matrice $M$ par le vacteur (colonne) $\ii$. D'autre part, si $M$ est dans $\GL_n(\Z)$, l'application $\x\mapsto ^M\x$ est une bijection de $(k^\times)^n$, ainsi que pour toutes ses extensions. Donc on obtient $L(^Mf;T)=L(f;T)$. 

Soit maintenant $f$ un polynôme de Laurent de la forme $f(\x)=\sum_{i=1}^n \sum_{j=-d_i'}^{d_i} a_{ij} \x^{j\e_i}$, dont le polyèdre convexe est $\Delta$. En choisissant pour $M$ la matrice de passage de la base $\{\e_1,\dots,\e_n\}$ \`a la base canonique, on voit que $^Mf(\x)= \sum_{i=1}^n \sum_{j=-d_i'}^{d_i} a_{ij} x_i^{j}$, dont le polyèdre associ\'e est l'enveloppe convexe des points de coordonnées $$(d_1,0,\dots,0),(-d_1',0,\dots,0),\dots,(0,\dots,0,d_n),(0,\dots,0,-d_n').$$ 
Mais ce dernier polyèdre est la somme directe des segments $[-d_i',d_i]$, $1\leq i\leq n$, et le corollaire \ref{prodgnp} nous assure que
$$\GNP([-d_1',d_1],p)\times \dots \times\GNP([-d_n',d_n],p)\preceq\GNP(\Delta,p).$$
D'autre part, d'apr\`es la proposition \ref{prodH}, on a $\HP(\Delta)=\HP([-d_1',d_1])\times \dots \times\HP([d_n',d_n])$. Le r\'esultat d\'ecoule maintenant du fait que pour chaque $i$, $\GNP([-d_i',d_i],p)$ tend vers $\HP([-d_i',d_i])$ quand $p$ tend vers $\infty$ : le polygone $\GNP(\Delta,p)$ est encadré entre deux polygones ayant la même limite.
\end{proof}
 
\begin{rem}
Dans le cas où $p\equiv 1$ modulo $\ppcm(d,d')$, on sait ({\it cf.} \cite{Ro}) que les polygones $\GNP([-d',d],p)$ et $\HP([-d',d])$ coïncident. En particulier, si $p\equiv 1$ modulo $D=\ppcm(d_i,d_i')_{1\leq i\leq n}$, on en déduit que les polygones $\GNP(\Delta,p)$ et $\HP(\Delta)$ coïncident, et la conjecture d'Adolphson et Sperber \cite[p. 386]{AS1} est vérifiée dans ce cas. 
\end{rem}
 
Considérons maintenant la seconde question, à savoir l'existence d'un ouvert dense $\U_{\Delta}$ défini sur $\Q$ de l'espace des polynômes à coefficients dans $\overline{\Q}$ de polyèdre de Newton $\Delta$ tel que pour tout $f$ de $\U_{\Delta}$, on ait $\lim_{p\rightarrow\infty} \NP_q(f\mod \p)=\HP(\Delta)$, où $\p$ est un premier au dessus de $p$ dans le corps $\Q_f$ engendré par les coefficients de $f$. Comme on ne considère pas tous les polynômes de poltope $\Delta$, on ne peut répondre à cette question. En revanche, pour les sous-familles que nous avons utilisées, on obtient le résultat suivant : 

\begin{theo}
\label{dense}
Il existe un ouvert  dense $\U$ défini sur $\Q$ de l'espace des polynômes de la forme $f(x)=\sum_{i=1}^n \sum_{j=-d_i'}^{d_i} A_{ij} \x^{j\textbf{e}_i}$ à coefficients dans $\overline{\Q}$ tel que pour tout polynôme dans $\U$, on ait 
$$\lim_{p\rightarrow\infty} \NP_q(f\mod \mathfrak{p})=\HP(\Delta).$$
\end{theo} 

\begin{proof}
Pour un polyn\^ome de la forme $f(x)=\sum_{i=1}^n \sum_{j=-d_i'}^{d_i} a_{ij} \x^{j\e_i}$, on a $\NP_q(f\mod \mathfrak{p})=\GNP([-d_1',d_1'],p)\times \dots \times\GNP([-d_n',d_n],p)$ si et seulement si les coefficients de $f$ v\'erifient
 $$\prod_{k=1}^n \H_{[-d_k',d_k]}^{\rho_k,\rho_k'}(A_{k,-d_k'},\dots, A_{k,d_k})\neq 0\mod p,$$
où $\rho_k$ ({\it resp.} $\rho_k'$) est le reste de la division euclidienne de $p$ par $d_k$ ({\it resp.} par $d_k'$). Notons $\H^{\rho}$ ce polyn\^ome, avec $\rho=(\rho_1,\rho_1',\dots,\rho_n,\rho_n')$. Comme on sait que les polyn\^omes $\H_{[-d_k',d_k]}^{\rho_k,\rho_k'}$ peuvent \^etre choisis \`a coefficients dans $\Q$, $\H^\rho$ est \`a coefficients dans $\Q$, et si 
$$\H(X_{ij})=\prod_{\rho \in (\z/d_1\z)^\times\times \dots\times (\z/d_n'\z)^\times}\H^{\rho}(X_{ij})=\prod_{i=1}^n \H_{[-d_i',d_i]}(X_{ij}),$$
on voit que pour tout $f$ dont les coefficients sont hors de l'hypersurface d'\'equation $\H=0$, on a pour tout $p$ assez grand 
$$\NP_q(f\mod \p)=\GNP([-d_1',d_1'],p)\times \dots \times\GNP([-d_n',d_n],p),$$ 
et le résultat découle de la convergence de ce dernier polygone vers $\HP(\Delta)$.
\end{proof}

\section{Comportement asymptotique, cas mixte.}
\label{HS}

Dans cette section, on note encore $(\e_1,\dots,\e_n)$ une base du $\Z$-module $\Z^n$, et $\Delta$ le polyèdre convexe de $\R^n$ qui est l'enveloppe convexe de l'ensemble des points $\{d_i\e_i,-d_i'\e_i\}_{1\leq i\leq n}$ et de l'origine si n\'ecessaire. On va étudier le comportement asymptotique de polygones de Newton de la forme $\NP_q(f,\chi)$ pour $f$ comme plus haut, et $\chi$ un caractère multiplicatif de $k^{\times}$ d'ordre fixé. Cette étude a été menée en dimension $1$ dans \cite{BFZ}, et nous allons la généraliser ici. Les résultats sont assez différents puisque qu'il n'y a plus de limite, mais une limite pour chaque classe inversible modulo l'ordre du caractère.

\medskip
 
Dans le cas où le caractère multiplicatif $\chi$ n'est pas trivial, la situation est assez différente. Soit $\omega$ le caractère de Teichmüller de $k^\times$, qui est un générateur de groupe des caractères de $k^\times$. Pour un $n$-uplet d'entiers $\dd=(\delta_1,\dots,\delta_n)$, on note $\chi=\omega^{\dd}$ le caractère de $(k^\times)^n$ défini par $\chi(x_1,\dots,x_n)=\omega(x_1)^{\delta_1}\dots \omega(x_n)^{\delta_n}$. Soit $f$ un polynôme de polyèdre de Newton $\Delta$ (qu'on suppose engendrer $\R^n$), non dégénéré. Adolphson et Sperber ont montré que la fonction $L(f,\chi;T)^{(-1)^{n-1}}$ est un polynôme de degré $n!V(\Delta)$ ; ils ont aussi donné une borne inférieure pour son polygone de Newton \cite[Theorem 3.17]{AS2}, qu'on appelera dans la suite {\it polygone de Hodge associé à $\Delta$ et $\dd$}, et qu'on notera $\HP(\Delta,\frac{\dd}{q-1})$. 

Nous allons décrire ce polygone, dans le cas où le polyèdre $\Delta$ engendre l'espace $\R^n$. Pour deux entiers $i$ et $0\leq\delta\leq q-2$, on note $\delta^{(i)}$ le reste modulo $q-1$ de l'entier $p^i\delta$ ; remarquons que la suite $(\delta^{(i)})_i$ est périodique : on a $\delta^{(a)}=\delta$, où $a=\log_p q$. On note encore $\dd^{(i)}=(\delta_1^{(i)},\dots,\delta_n^{(i)})$.

Soit maintenant $N^{(i)}$ le réseau $\frac{\dd^{(i)}}{q-1}+\Z^n$ de $\R^n$. On note $M_{\Delta,\dd^{(i)}}:=C(\Delta)\cap N^{(i)}$, et $\AA_{\Delta,\dd^{(i)}}$ le $\AA_{\Delta}$-module $k[x^{M_{\Delta,\dd^{(i)}}}]$. Il existe alors un entier positif $D$, minimal, tel que l'image de $M_{\Delta,\dd^{(i)}}$ par $w_{\Delta}$ soit contenue dans $\frac{1}{D}\N$ ; on appellera cet entier le {\it dénominateur de $(\Delta,\dd^{(i)})$} dans la suite. Muni de ce poids, $\AA_{\Delta,\dd^{(i)}}$ devient un $\AA_{\Delta}$-module gradué, auquel on peut associer une série de Poincaré et un polynôme $P_{\Delta,\dd^{(i)}}$ comme plus haut. Notons $\Pi^{(i)}$ le polygone issu de cette série. De plus chacun des polynômes $P_{\Delta,\dd^{(i)}}$ est de degré plus petit que $nD$, et vérifie $P_{\Delta,\dd^{(i)}}(1)=n!V(\Delta)$. On a ainsi une famille de polygones $\Pi^{(i)}$ pour $0\leq i\leq a$, tous de même longueur $n!V(\Delta)$. 

Si $\Pi$ et $\Pi'$ sont deux polygones de même longueur, on note $\Pi+\Pi'$ le polygone dont la pente sur le segment $[i,i+1]$ est la somme des pentes de $\Pi$ et $\Pi'$ sur ce segment ; d'autre part, pour un réel $r>0$, on désigne par $r\Pi$ le polygone obtenu à partir de $\Pi$ en multipliant toute ses pentes par $r$ (c'est en fait l'image de $\Pi$ par l'affinité orthogonale d'axe $Ox$, de direction $Oy$ et de rapport $r$). Avec ces notations, on sait alors décrire le polygone de Hodge
$$ \HP(\Delta,\frac{\dd}{q-1})=\frac{1}{a}\sum_{i=0}^{a-1}\Pi^{(i)}~;$$

On sait maintenant exprimer le polygone de Hodge associé à une somme directe de polyèdres et à deux caractères multiplicatifs, à l'aide des polygones de Hodge associés à ses facteurs. C'est l'exacte transposition de la proposition \ref{prodH} dans ce nouveau cadre, et nous en omettons la démonstration, qui est très similaire au cas du polygone de Hodge de l'algèbre $\AA_{\Delta}$.

\begin{propo}
\label{prodHS}
Soient $\Delta_1$ et $\Delta_2$ deux polyèdres convexes, respectivement dans $\R^{n_1}$ et $\R^{n_2}$, et $\Delta$ leur somme directe. Si on pose $$\dd_1=(\delta_1,\dots,\delta_{n_1}),~\dd_2=(\delta_{n_1+1},\dots,\delta_{n_1+n_2}),~\mbox{\it et}~\dd=(\dd_1,\dd_2)=(\delta_1,\dots,\delta_{n_1+n_2}),$$ 
alors le polygone de Hodge $\HP(\Delta,\frac{\dd}{q-1})$ est le produit des polygones de Hodge de ses facteurs
$$\HP(\Delta,\frac{\dd}{q-1})=\HP(\Delta_1,\frac{\dd_1}{q-1})\times \HP(\Delta_2,\frac{\dd_2}{q-1}).$$
\end{propo}

Le cas de la dimension $1$ a été étudié dans \cite{BFZ} ; le lecteur intéressé par les preuves des résultats que nous allons rappeler maintenant pourra s'y référer. Dans ce cas, le polygone de Hodge fait intervenir à la fois le polygone de Hodge des sommes associées à un caractère additif, et la valuation de la somme de Gauss associée à $\chi$, donnée par le théorème de Stickelberger ; pour cette raison nous avons choisi de renommer ce polygone {\it polygone de Hodge Stickelberger}, en changeant légèrement les notations. 

Pour justifier ces changements, commençons par expliquer le point de vue de \cite{BFZ} plus précisément. Les auterus sont motivés par l'étude du comportement asymptotique dans un cas non générique des polynômes de Laurent en une variable (plus particulièrement par le cas des polynômes de la forme $f(x^s)$). La formule de Poisson ramène ce problème à la situation suivante : le polynôme de Laurent $f\in k[x,x^{-1}]$ admet pour polyèdre de Newton le segment $[-d',d]$, et $\chi$ est un caractère multiplicatif d'ordre divisant $s$. Ce dernier n'est pas nécessairement défini sur $k$, mais sur une de ses extensions $k'$ de cardinal $q'\equiv 1$ modulo $s$. Si $\omega'$ est le caractère de Teichmüller de $k'^\times$, on peut donc écrire $\chi=\omega'^{\delta}$, avec $\delta=\frac{(q'-1)r}{s}$ pour un certain entier $1\leq r\leq s-1$. Alors le polygone de Hodge $\HP([-d',d],\frac{\delta}{q'-1})$ est le polygone d'extr\'emit\'es l'origine et le point de coordonn\'ees $(d+d',\frac{d+d'}{2})$, et poss\'edant un segment de longueur $1$ pour chacune des pentes suivantes
$$\frac{1-\lambda}{d},\dots,\frac{d-\lambda}{d},\frac{\lambda}{d'},\dots,\frac{d'-1+\lambda}{d'} ~\left( \frac{1-\lambda}{d},\dots,\frac{d-\lambda}{d} ~\mbox{\rm si}~d'=0\right),$$
où, comme dans le théorème de Stickelberger donnant la valuation des sommes de Gauss, $\lambda=\frac{1}{\log_p(q')(p-1)}s_p\left(\frac{(q'-1)r}{s}\right)$ avec $s_p$ la somme des chiffres de l'écriture en base $p$ de l'entier $\frac{(q'-1)r}{s}$. En particulier, $\lambda$ ne dépend pas du choix de $q'$. 

On sait réexprimer $\lambda$ de la façon suivante : si $\sigma_p$ désigne la permutation de l'ensemble $\{0,\dots,s-1\}$ induite par la multiplication par $p$ dans $\Z/s\Z$, et si $\sigma$ est le cycle, de longueur $\ell$, de $\sigma_p$ contenant $r$, alors on a $\lambda=\frac{\sum_{j\in \sigma} j}{s\ell}$. Donc $\lambda$ ne dépend pas de $p$, mais de son reste modulo $s$. 

\begin{defi}
Le polygone qu'on vient de décrire est le polygone de Hodge Stickelberger associé au polyèdre $[-d',d]$ et au rationnel $\frac{r}{s}$. On le note $\HS([-d',d],\frac{r}{s},\nu)$, où $\nu$ est le reste de $p$ modulo $s$.
\end{defi}
  
Ici encore, on a un polygone de Newton générique $\GNP([-d',d],\chi,p)$, et un polynôme de Hasse à coefficients rationnels ; ils sont eux aussi indépendants de la puissance de $p$ choisie. De plus, le polynôme de Hasse ne dépend que des restes respectifs $\nu,\rho$ et $\rho'$ de $p$ modulo $s$, $d$ et $d'$ ; notons le $\H_{[-d',d],\frac{r}{s},\nu}^{\rho,\rho'}$. 

On voit que la principale différence avec les sommes additives est que, pour un polyèdre fixé, il y a plusieurs polygones de Hodge-Stickelberger, dépendant du reste de la division euclidienne de $p$ par l'ordre du caractère multiplicatif $\chi$. 

En conséquence, on ne peut plus espérer que les polynômes de Newton génériques convergent quand $p$ tend vers $+\infty$. En revanche, quand $p$ tend vers $+\infty$ le ``long d'une classe de $(\Z/s\Z)^\times$", on obtient la convergence :
$$\lim_{\stackrel{p\rightarrow +\infty}{p\equiv\nu~[s]}} \GNP([-d',d],\chi,p)=\HS([-d',d],\frac{r}{s},\nu).$$
D'autre part, en définissant  $\H_{[-d',d],\frac{r}{s},\nu}=\prod_{\rho,\rho'}\H_{[-d',d],\frac{r}{s},\nu}^{\rho,\rho'}$, on sait que pour tout polynôme $f\in\overline{\Q}[x,x^{-1}]$, de la forme $f(x)=\sum_{i=-d'}^d A_ix^i$, et vérifiant 
$$\H_{[-d',d],\frac{r}{s},\nu}(A_{-d'},\dots,A_d)\neq 0,$$ 
on a $\lim_{\stackrel{p\rightarrow +\infty}{p\equiv\nu~[s]}} \NP_q(f\mod \p,\chi)=\HS([-d',d],\frac{r}{s},\nu)$ (où comme plus haut $\p$ est un premier au dessus de $p$ dans le corps $\Q_f$ engendré par les coefficients de $f$). C'est à dire qu'il existe un ouvert dense $\U_{\frac{r}{s},\nu}$, défini sur $\Q$, de l'espace des polynômes de polyèdre $[-d',d]$ en une variable à coefficients dans $\overline{\Q}$, tel que pour tout $f$ dans $U_{\frac{r}{s},\nu}$, la limite $\lim_{\stackrel{p\rightarrow +\infty}{p\equiv\nu~[s]}} \NP_q(f\mod \p,\chi)$ existe et coïncide avec le polygone de Hodge-Stickelberger.

\medskip

On peut maintenant passer de la dimension $1$ à la dimension supérieure. Commençons par définir les polygones de Hodge-Stickelberger dans ce cadre :

\begin{defi}
On pose, pour $\Delta$ un polyèdre convexe, $\frac{\r}{\s}=(\frac{r_1}{s_1},\dots,\frac{r_n}{s_n})$, $p$ un premier de résidu $\nu$ modulo $s=\ppcm(s_1,\dots,s_n)$ et $q$ une puissance de $p$ telle que $q\equiv 1\mod s$
$$\HS(\Delta,\frac{\r}{\s},\nu):=\HP(\Delta,\frac{\dd}{q-1})$$
avec $\dd=\left(\frac{(q-1)r_1}{s_1},\dots,\frac{(q-1)r_n}{s_n}\right)$.
\end{defi}

\begin{rem}
On vérifie aisément, à l'aide de la définition qu'on en a donné, que ce polygone ne dépend pas de $q$, la puissance de $p$ qu'on choisit. D'autre part il ne dépend ici encore que du reste de $p$ modulo $s$, ce qui justifie notre notation.
\end{rem}

La proposition \ref{prodHS} s'applique de la façon suivante aux polygones de Hodge-Stickleberger ; si 
\begin{itemize}
	\item[\it i/] $\Delta_1$ et $\Delta_2$ sont deux polèdres convexes, et $\Delta=\Delta_1\oplus\Delta_2$ leur somme directe ;
	\item[\it ii/] $\frac{\r_1}{\s_1}$, $\frac{\r_2}{\s_2}$ et $\frac{\r}{\s}=\left(\frac{\r_1}{\s_1},\frac{\r_2}{\s_2}\right)$ des $n_1$, $n_2$ et $n$-uplets de rationnels ;
	\item[\it iii/] $\nu$ un résidu inversible modulo $s:=\ppcm(s_1,s_2)$, et $\nu_i$ son image modulo $s_i$, alors
\end{itemize}
  $$\HS(\Delta,\frac{\r}{\s},\nu)=\HS(\Delta_1,\frac{\r_1}{\s_1},\nu_1)\times \HS(\Delta_2,\frac{\r_2}{\s_2},\nu_2).$$

On rappelle que $(\e_1,\dots,\e_n)$ est une base du $\Z$-module $\Z^n$, et que $\Delta$ désigne le polyèdre convexe de $\R^n$ qui est l'enveloppe convexe de l'ensemble des points $\{d_i\e_i,-d_i'\e_i\}_{1\leq i\leq n}$ et de l'origine si n\'ecessaire. D'autre part, on choisit, pour chaque premier $p$ assez grand, un caractère multiplicatif $\chi_i=\omega_{q-1}^{\frac{(q-1)r_i}{s_i}}$ pour $1\leq i\leq n$, avec $q$ une puissance convenable de $p$. On note $\chi=(\chi_1,\dots,\chi_n)$ le caractère de $(\F_q^\times)^n$ induit par les $\chi_i$, et $\frac{\r}{\s}=(\frac{r_1}{s_1},\dots,\frac{r_n}{s_n})$. Alors le polygone 
$$\HS(\Delta,\frac{\r}{\s},\nu):=\HS([-d_1',d_1],\frac{r_1}{s_1},\nu_1)\times \dots \times \HS([-d_n',d_n],\frac{r_n}{s_n},\nu_n)$$
ne dépend que du reste de la division euclidienne de $p$ par $s:=\ppcm(s_1,\dots,s_n)$. A l'aide de ces notations, les transpositions des théorèmes \ref{lim} et \ref{dense} à cette nouvelle situation s'écrivent :

\begin{theo}
\label{limS}
Quand $p$ tend vers l'infini, le polygone de Newton g\'en\'erique de $\Delta$ associ\'e au premier $p$ et au caractère $\chi$, $\GNP(\Delta,\chi,p)$, tend vers le polygone de Hodge-Stickelberger $\HS(\Delta,\frac{\r}{\s},\nu)$ quand $p$ tend vers $+\infty$ avec $p\equiv \nu~[s]$.
\end{theo}

\begin{theo}
\label{denseS}
Il existe un ouvert  dense $\U$ défini sur $\Q$ de l'espace des polynômes de la forme $f(x)=\sum_{i=1}^n \sum_{j=-d_i'}^{d_i} a_{ij} \x^{j\e_i}$ à coefficients dans $\overline{\Q}$ tel que pour tout polynôme dans $\U$, on ait 
$$\lim_{\stackrel{p\rightarrow +\infty}{p\equiv\nu~[s]}} \NP_q(f\mod \mathfrak{p},\chi_p)=\HS(\Delta,\frac{r}{s},\nu).$$
\end{theo} 

\begin{rem}
Dans le cas où $p\equiv 1$ modulo $\ppcm(d,d',s)$, on sait ({\it cf.} \cite{BFZ}) que les polygones $\GNP([-d',d],\chi,p)$ et $\HS([-d',d],\frac{r}{s},1)$ coïncident. En particulier, si $p\equiv 1$ modulo $D=\ppcm(d_i,d_i',s_i)_{1\leq i\leq n}$, on en déduit que les polygones $\GNP(\Delta,\chi,p)$ et $\HP(\Delta,\frac{\r}{\s},1)$ coïncident, et on a vérifié dans ce cas une extension de la conjecture d'Adolphson et Sperber aux sommes mixtes. 
\end{rem}

Nous terminons ce chapitre par le cas particulier $s=2$. Pour un nombre premier impair $p$ fixé, on note $\chi_2$ le caractère quadratique, défini sur $\F_q^\times$ par $\chi_2(x)=\omega^{\frac{q-1}{2}}(x)$. Tous les caractères multiplicatifs d'ordre $2$ de $(\F_q^\times)^n$ sont de la forme $\chi_2^{\ee}$, avec $\ee=(\ep_1,\dots,\ep_n)\in \{0,1\}^n$, et $\chi_2^{\ee}(x_1,\dots,x_n)=\chi_2^{\ep_1}(x_1)\dots\chi_2^{\ep_n}(x_n)$. Puisque (presque) tous les premiers sont impairs, le polygone de Hodge-Stickelberger ne dépend plus que de $\ee$, on le note $\HS(\Delta,\frac{\ee}{2})$. Pour la même raison, on peut décrire directement ce polygone à l'aide d'une série de Poincaré. 

\begin{lem}
\label{ee2}
Soit $\Delta$ un polèdre convexe de $\R^n$ qui l'engendre, $\AA_{\Delta,\frac{\ee}{2}}$ le $\AA_{\Delta}$-module gradué associé à cette situation. Alors le polygone $\HS(\Delta,\frac{\ee}{2})$ est le polygone issu de la série de Poincaré de $\AA_{\Delta,\frac{\ee}{2}}$.
\end{lem}

L'indépendance du polygone de Hodge par rapport à $p$ nous permet de retrouver l'existence d'une limite. 

\begin{coro}
\label{limS2}
Quand $p$ tend vers l'infini, le polygone de Newton g\'en\'erique de $\Delta$ associ\'e au premier $p$ et au caractère quadratique $\chi_2$, $\GNP(\Delta,\chi_2,p)$, tend vers le polygone de Hodge-Stickelberger $\HS(\Delta,\frac{\ee}{2})$ quand $p$ tend vers $+\infty$.
\end{coro}

\begin{coro}
\label{denseS2}
Il existe un ouvert  dense $\U$ défini sur $\Q$ de l'espace des polynômes de la forme $f(x)=\sum_{i=1}^n \sum_{j=-d_i'}^{d_i} a_{ij} \x^{j\e_i}$ à coefficients dans $\overline{\Q}$ tel que pour tout polynôme dans $\U$, on ait 
$$\lim_{p\rightarrow +\infty} \NP_q(f\mod \mathfrak{p},\chi_2)=\HS(\Delta,\frac{\ee}{2}).$$
\end{coro}

\section{Polynômes de polyèdres d'exposant deux.}

Dans cette partie on va étendre les principaux résultats à des polyèdres un peu plus généraux : on fixe un entier $n$, et on note $(\f_1,\dots,\f_n)$ une famille libre de $\Z^n$, qui engendre un sous module $N$ tel que le quotient $\Z^n/N$ soit un groupe d'exposant $2$. Comme dans les chapitres précédents on choisit des entiers naturels $d_1,d_1',\dots,d_n,d_n'$. On note $\Delta$ le polyèdre convexe de $\R^n$ qui est l'enveloppe convexe de l'ensemble des points $\{d_i\f_i,-d_i'\f_i\}_{1\leq i\leq n}$ et de l'origine si n\'ecessaire. On va réexprimer les sommes additives associées à certains polynômes de polyèdre $\Delta$ à l'aide de sommes mixtes étudiées dans le chapitre précédent ; on utilisera ensuite les corollaires \ref{limS2} et \ref{denseS2} pour obtenir la limite.

On se propose de d\'emontrer les r\'esultats suivants :

\begin{theo}
\label{lim2}
Quand $p$ tend vers l'infini, le polygone de Newton g\'en\'erique de $\Delta$ associ\'e au premier $p$, $\GNP(\Delta,p)$, tend vers le polygone de Hodge $\HP(\Delta)$.
\end{theo}

\begin{theo}
\label{dense2}
Il existe un ouvert  dense $\U$ défini sur $\Q$ de l'espace des polynômes de la forme $f(x)=\sum_{i=1}^n \sum_{j=-d_i'}^{d_i} a_{ij} \x^{j\textbf{f}_i}$ à coefficients dans $\overline{\Q}$ tel que pour tout polynôme dans $\U$, on ait 
$$\lim_{p\rightarrow\infty} \NP_q(f\mod \mathfrak{p})=\HP(\Delta).$$
\end{theo} 

\medskip

Dans toute la suite, on suppose que $p$ est un nombre premier impair.

\medskip

Soit $\FF=(\f_1,\dots,\f_n)$ une famille libre de $\Z^n$, qui engendre un sous module $N$ de $\Z^n$ tel que le quotient $\Z^n/N$ soit un groupe d'exposant $2$. On note $M=(f_{ij})$ la matrice de passage de la base canonique à la famille $\FF$ dans $\MM_n(\Z)$, et $k$ la dimension comme $\F_2$-espace vectoriel de $\Z^n/N$. On a donc la suite exacte
\begin{equation}
\label{quot}
0\rightarrow \Z^n \rightarrow \Z^n\rightarrow \Z^n/N \simeq \F_2^k\rightarrow 0,
\end{equation}
où la première flèche est l'action de $M$.

On peut trouver une base $(\e_1,\dots,\e_n)$ de $\Z^n$ (comme $\Z$-module) telle que la famille $\e_1,\dots,\e_{n-k},2\e_{n-k+1},\dots,2\e_n$ soit une base de $N$. De façon équivalente, la matrice $M$ est équivalente, sur $\MM_n(\Z)$, à la matrice diagonale dont les $n-k$ premiers coefficients diagonaux valent $1$, et les $k$ derniers valent $2$ ; remarquons qu'on doit avoir $\det M=2^k$.

\begin{defi}
On note $M_2$ l'application linéaire de $\F_2^n$ induite par la réduction modulo $2$ de la matrice $M$. Soit $\overline{E}$ son noyau, et 
pour tout $\overline{\ee}=(\overline{\ep}_1,\dots,\overline{\ep}_n)\in \overline{E}$, on note $\ee=(\ep_1,\dots,\ep_n)$ le relèvement de $\overline{\ee}$ à $\{0,1\}^n$. Finalement, soit $E$ le sous ensemble de $\{0,1\}^n$ formé des $\ee$ quand $\overline{\ee}$ décrit $\overline{E}$.
\end{defi}

L'ensemble $E$ qu'on vient de définir va nous servir à décrire les points entiers d'un domaine fondamental de $\Z^n/N$. D'autre part, rappelons que si $\x=(x_1,\dots,x_n)$, on note $^M\x=(y_1,\dots,y_n)$, avec $y_i=\prod_{j=1}^n x_j^{f_{ji}}$. L'ensemble $E$ va aussi nous servir à décrire l'image du morphisme $\varphi_M:\x\mapsto ^M\x$ de $k^{\times n}$ dans lui-même, et à réexprimer les sommes pures, (ainsi que leurs fonctions $L$ et groupes de cohomologie) associées au polynôme $f(^M\x)$ en fonction de sommes mixtes associèes au polynôme $f$ et à certains caractères quadratiques. On rappelle qu'on note $\chi_2^{\ee}$, avec $\ee=(\ep_1,\dots,\ep_n)\in \{0,1\}^n$, le caractère (multiplicatif) de $k^{\times n}$ défini par $\chi_2^{\ee}(x_1,\dots,x_n)=\chi_2^{\ep_1}(x_1)\dots\chi_2^{\ep_n}(x_n)$.

\begin{lem}
\label{E}
\begin{itemize}
	\item[ i/] L'ensemble des points entiers contenus dans le polyèdre 
	$$[0,1[\f_1\times \dots\times [0,1[\f_n=\left\{\sum_{i=1^n}x_i\f_i,~0\leq x_i<1\right\}$$
	est $\{ f_{\ee}:=\frac{1}{2}\sum_{i=1}^n \ep_i \f_i,~\ee\in E\}$.
	\item[ii/] Le sous-groupe du groupe des caractères multiplicatifs de $k^{\times n}$ orthogonal à l'image du morphisme $\varphi_M$ est 
	$$(\Ima \varphi_M)^{\perp}=\left\{\chi_2^{\ee},~\ee\in E\right\}.$$
\end{itemize}
\end{lem}

\begin{proof}
Le polyèdre $[0,1[\f_1\times \dots\times [0,1[\f_n$ est un domaine fondamental pour l'action de $N$ sur $\R^n$ par translations. En particulier il contient $\det M=2^k$ points entiers. Puisque $M_2$ est de rang $k$, $E$ contient $2^k$ éléments, et il suffit de vérifier que les points $\frac{1}{2}\sum_{i=1}^n \ep_i \f_i$ sont entiers. Mais par construction l'image de $\sum_{i=1}^n \ep_i \f_i$ dans $\F_2^n$ est nulle, c'est à dire que toutes les coordonnées de $\sum_{i=1}^n \ep_i \f_i$ sont paires, et on a prouvé l'assertion {\it i/}.

On a supposé $p$ impair, et le groupe $k^{\times n}$ est isomorphe à $\left(\Z/(q-1)\Z\right)^n$, avec $q-1$ pair. En tensorisant la suite exacte (\ref{quot}) par le groupe $\Z/(q-1)\Z$, on obtient donc la suite exacte
\begin{equation}
\label{quot2}
\left(\Z/(q-1)\Z\right)^n \rightarrow \left(\Z/(q-1)\Z\right)^n \rightarrow  \F_2^k \rightarrow 0.
\end{equation}
Ainsi l'image du morphisme $\varphi_M$ est d'indice $2^k$ dans $k^{\times n}$, et il suffit encore une fois de montrer une inclusion. Un calcul facile montre que $\chi_2^{\ee}(^M\x)=\chi_2^{e_1}(x_1)\dots\chi_2^{e_n}(x_n)$, où $e_i$ est la $i$-ème coordonnée du vecteur $M\ee$. Mais comme plus haut, quand $\ee$ est dans $E$, le vecteur $M\ee$ a toutes ses coordonnées paires. Donc $\chi_2^{\ee}$ est dans l'orthogonal de l'image de $\varphi_M$, ce qui termine la démonstration.
\end{proof}

Nous pouvons maintenant réexprimer les sommes pures associées à un polynôme de Laurent de la forme $f(^M\x)$.

\begin{propo}
\label{poisson}
Soit $f\in k[\x,\x^{-1}]$ un polynôme de Laurent, et $M$ comme précédemment ; on rappelle que $^Mf$ est le polynôme de Laurent $f(^Mx)$. On a les décompositions suivantes
\begin{itemize}
	\item[i/] des sommes de caractères
	$$\sum_{\textbf{x} \in k_r^{\times n}} \psi(^Mf(\x))=\sum_{\ee\in E}\sum_{\textbf{x} \in k_r^{\times n}} \psi(f(\x))\chi_2^{\ee}(\x)~;$$
	\item[ii/] de la fonction $L$ 
	$$L(^Mf,T)=\prod_{\ee\in E} L(f,\chi_2^{\ee},T)~;$$
	\item[iii/] de la suite exacte longue de cohomologie 
	$$H^\bullet_c(\G_m^n,(^Mf)^*\L_\psi)=\bigoplus_{\ee\in E} H^\bullet_c(\G_m^{n},f^*\L_\psi\otimes \L_{\chi_2^{\ee}}).$$
\end{itemize}
\end{propo}

\begin{proof}
Ce sont différents avatars de la formule de Poisson appliquée à notre situation. Montrons le premier : à l'aide de l'assertion {\it ii/} du lemme \ref{E}, on a 
$$\begin{array}{rcl}
\sum_{\textbf{x} \in k^{\times n}} \psi(f(^M\x)) & = & \# \Ker \varphi_M\sum_{\textbf{y}\in \textrm{Im}\varphi_M} \psi(f(\y)) \\
 & = & \# \Ker \varphi_M\left(\frac{1}{\#(\textrm{Im}\varphi_M)^{\perp}} \sum_{\chi\in (\textrm{Im}\varphi_M)^{\perp}}\sum_{\textbf{x}\in k^{\times n}} \psi(f(\x))\chi_2^{\ee}(\x)\right) \\
 & = & \sum_{\ee\in E}\sum_{\textbf{x} \in k_r^{\times n}} \psi(f(\x))\chi_2^{\ee}(\x),\\
 \end{array}$$
ce qui est exactement la formule recherchée.
\end{proof}

On va en déduire une décomposition pour le polygone de Newton ; commençons par introduire une nouvelle opération sur les polygones :

\begin{defi}
Si $\Pi_1$ et $\Pi_2$ sont deux polygônes convexes de pentes respectives $(s_i)_{1\leq i\leq a}$ et $(s_i')_{1\leq i\leq b}$, leur juxtaposition est le polygone convexe $\Pi_1\coprod\Pi_2$ de pentes $(s_i,s_j')_{1\leq i\leq a,1\leq j\leq b}$.
\end{defi}

A l'aide de cette définition, on déduit de la proposition \ref{poisson} l'écriture suivante :

\begin{coro}
\label{dec}
Le polygone de Newton associé au polynôme $^Mf$ est la juxtaposition des polygônes de Newton associés au polynôme $f$ et aux caractères $\chi_2^{\ee}$ quand $\ee$ décrit $E$
$$\NP_q(^Mf)=\coprod_{\ee\in E} \NP_q(f,\chi_2^{\ee}).$$
\end{coro}

Nous allons maintenant réécrire le polygone de Hodge $\HP(\Delta)$ à l'aide des polygones des sections précédentes :

\begin{lem}
\label{decH}
Soit $\Delta_0$ le polyèdre convexe de $\R^n$ qui est l'enveloppe convexe de l'ensemble des points $\{d_i\e_i,-d_i'\e_i\}_{1\leq i\leq n}$ et de l'origine si n\'ecessaire, pour $(\e_1,\dots,\e_n)$ une base du $\Z$ module $\Z^n$. On a la décomposition suivante du polygone de Hodge associé à $\Delta$ :
$$\HP(\Delta)=\HP(\Delta_0)\coprod\left(\coprod_{\ee\in E\backslash\{0,\dots,0\}} \HS(\Delta_0,\frac{\ee}{2})\right).$$
\end{lem}

\begin{proof}
On va revenir ici à la définition du polygône de Hodge à l'aide des séries de Poincaré des algèbres $\AA_{\Delta}$ et $\AA_{\Delta_0,\frac{\ee}{2}}$.
Quitte à permuter les couples $(d_i,d_i')$, et à échanger $d_i$ et $d_i'$ dans certains couples, on peut supposer que $d_1'=\dots=d_l'=0$ et que les $d_i'$ sont tous non nuls pour $l+1\leq i\leq n$. Alors d'après le lemme \ref{E}, les points de $M_{\Delta}$ sont les 
$$\f_{\textbf{k},\ee}:=k_1\f_1+\dots+k_n\f_n+\f_{\ee},~\k=(k_1,\dots,k_n)\in \N^l\times\Z^{n-l},~\ee\in E,$$
et le poids d'un tel point est donné par
$$w_{\Delta}(\f_{\textbf{k},\ee})=\sum_{i=1}^l\frac{k_i+\frac{\ep_i}{2}}{d_i}+\sum_{i=l+1}^{n}\max\left(\frac{k_i+\frac{\ep_i}{2}}{d_i},-\frac{k_i+\frac{\ep_i}{2}}{d_i'}\right).$$
On en déduit, en notant $D$ le dénominateur de $\Delta$, que la série de Poincaré de l'algèbre $\AA_{\Delta}$ s'écrit :
$$P_{\AA_{\Delta}}(t)=\sum_{\ee\in E}\sum_{\textbf{k}\in \n^l\times\z^{n-l}}t^{Dw_{\Delta}(\textbf{f}_{\textbf{k},\ee})}.$$
Fixons $\ee\in E$. Les points de $M_{\Delta_0,\frac{\ee}{2}}$ sont les $\e_{\textbf{k},\ee}=k_1\e_1+\dots+k_n\e_n+\e_{\ee}$, où $\k$ décrit $\N^l\times\Z^{n-l}$ et $\e_{\ee}:=\frac{1}{2}\sum_{i=1}^n\ep_i\e_i$, et on a 
$$w_{\Delta_0}(\e_{\textbf{k},\ee})=\sum_{i=1}^l\frac{k_i+\frac{\ep_i}{2}}{d_i}+\sum_{i=l+1}^{n}\max\left(\frac{k_i+\frac{\ep_i}{2}}{d_i},-\frac{k_i+\frac{\ep_i}{2}}{d_i'}\right).$$
Si $D_{\ee}$ désigne le dénominateur de $(\Delta_0,\frac{\ee}{2})$, la série de Poincaré de $\AA_{\Delta_0,\frac{\ee}{2}}$ s'écrit donc
$$P_{\AA_{\Delta_0,\frac{\ee}{2}}}(t)=\sum_{\textbf{k}\in \n^l\times\z^{n-l}}t^{D_{\ee}w_{\Delta}(\textbf{f}_{\textbf{k},\ee})}.$$
Remarquons au passage que $w_{\Delta_0}(M_{\Delta_0,\frac{\ee}{2}})\subset w_{\Delta}(M_{\Delta})$, c'est à dire que $D_{\ee}$ divise $D$. On peut donc écrire :
$$P_{\AA_{\Delta}}(t)=\sum_{\ee\in E} P_{\AA_{\Delta_0,\frac{\ee}{2}}}(t^{\frac{D}{D_{\ee}}}),$$
et en multipliant les deux membres par $(1-t^D)^n$, on trouve
$$P_{\Delta}(t)= P_{\Delta_0}(t^{\frac{D}{D_0}})+\sum_{\ee\in E\backslash\{0,\dots,0\}} P_{\Delta_0,\frac{\ee}{2}}(t^{\frac{D}{D_{\ee}}}).$$

Le résultat découle maintenant de la construction du polygone de Hodge de $\Delta_0$, et du lemme \ref{ee2} qui relie les polygones $\HS(\Delta_0,\frac{\ee}{2})$ aux séries de Poincaré des $\AA_{\Delta_0,\frac{\ee}{2}}$.
\end{proof}

Ces résultats intermédiaires nous donnent tous les éléments pour démontrer les théorèmes \ref{lim2} et \ref{dense2}.

\begin{proof}{(du théorème \ref{lim2})}
Remarquons d'abord que 
$$f(\x)=\sum_{i=1}^n \sum_{j=-d_i'}^{d_i} a_{ij} \x^{j\textbf{f}_i}=^Mg(x),$$
pour le polynôme $g(\x)=\sum_{i=1}^n \sum_{j=-d_i'}^{d_i} a_{ij} \x_i^{j}$. D'après le corollaire \ref{dec}, le polygone de Newton générique de la famille des polynômes $f(\x)=\sum_{i=1}^n \sum_{j=-d_i'}^{d_i} a_{ij} \x^{j\textbf{f}_i},~a_{ij}\in k$ est donné par 
$$\GNP(\Delta_0,p)\coprod\left(\coprod_{\ee\in E\backslash\{0,\dots,0\}} \GNP(\Delta_0,\frac{\ee}{2},p)\right).$$
Le théorème de spécialisation de Grothendieck assure que le polygone de Newton générique de la famille de tous les polynômes de polyèdre $\Delta$ vérifie
$$\GNP(\Delta_0,p)\coprod\left(\coprod_{\ee\in E\backslash\{0,\dots,0\}} \GNP(\Delta_0,\frac{\ee}{2},p)\right)\preceq\GNP(\Delta,p)\preceq \HP(\Delta).$$
Finalement le lemme \ref{decH}, joint au théorème \ref{lim} appliqué au polyèdre $\Delta_0$, et au corollaire \ref{limS2} nous assurent que le membre de gauche tend vers le membre de droite quand $p$ tend vers $\infty$. Donc le polygone $GNP(\Delta,p)$ tend vers $\HP(\Delta)$ quand $p$ tend vers $\infty$. 
\end{proof}

\begin{proof}{(du théorème \ref{dense2})}
Pour $f$ un polynôme de la forme 
$$f(\x)=\sum_{i=1}^n \sum_{j=-d_i'}^{d_i} a_{ij} \x^{j\textbf{f}_i}$$ 
à coefficients dans $\overline{\Q}$, on vérifie comme plus haut que
$$\lim_{p\rightarrow\infty} \NP_q(f\mod \mathfrak{p})=\HP(\Delta)$$
dès que $\H(a_{ij})\neq 0$, avec le polynôme de Hasse
$$\H(X_{ij})=\prod_{i=1}^n\prod_{\ee\in E} \H_{[-d_i',d_i],\frac{\ep_i}{2}}(X_{ij}).$$
\end{proof}


\begin{thebibliography}{99}


\bibitem{AS1}
{\sc A. Adolphson, S. Sperber :} Exponential sums and Newton polyhedra: Cohomology and estimates, 
{\it Ann. Math} {\bf 130} (1989), 367--408.

\bibitem{AS3}
{\sc A. Adolphson, S. Sperber:} On twisted exponential sums,
{\it Math. Ann.} {\bf 290} (1991), 713-726.

\bibitem{AS2}
{\sc A. Adolphson, S. Sperber :} Twisted exponential sums and Newton polyhedra, 
{\it J. reine angew. Math.} {\bf 443} (1993), 151--177.

\bibitem{BF}
{\sc R. Blache, \'E. F\'erard :}
Newton stratification for polynomials: the open stratum,
{\it J. Number Th.} {\bf 123} (2007), 456--472.

\bibitem{BFZ}
{\sc R. Blache, \'E. F\'erard, H.J.Zhu :}
Hodge-Stickelberger polygons for L-functions of exponential sums of $P(x^s)$,
preprint, available at  {\tt http://arxiv.org/abs/0706.2340}.

\bibitem{De}
{\sc P. Deligne :}
La conjecture de Weil : I.
{\it Publ. Math. IHES} {\bf 43} (1974), 273--307.

\bibitem{DL}
{\sc J. Denef, F. Loeser :}
Weights of exponential sums, intersection cohomology, and Newton polyhedra.
{\it Inv. Math.} {\bf 106} (1991), 275--294.

\bibitem{Dw}
{\sc B. Dwork:}
On the zeta function of a hypersurface.
{\it Publ. Math. IHES}{\bf 12}
(1962), 5--68.

\bibitem{Gr}
{\sc A. Grothendieck :} Formule de Lefschetz et rationalité des fonctions $L$, Séminaire Bourbaki, exposé 279, 1964/65.

\bibitem{hrz}
{\sc M. Henk, J. Richter-Gebert, G. Ziegler :} Basic properties of convex polytopes, dans Handbook of Discrete and Computational Geometry, CRC Press, 1997.

\bibitem{Ka}
{\sc N.M. Katz :} Slope filtration of $F$-crystals, {\it Ast\'erisque} {\bf 63} (1979), 113-164.

\bibitem{Kob}
{\sc N. Koblitz:}
$p$-adic numbers, $p$-adic analysis, and zeta-functions,
(Second edition),
{\it Graduate Texts in Mathematics} {\bf 58}.
Springer-Verlag, New York, 1984.

\bibitem{Kou}
{\sc A.G. Kouchnirenko :}
Poly\`edres de Newton et nombres de Milnor.
{\it Inv. Math.} {\bf 32} (1976), 1--31.

\bibitem{LZ}
{\sc H. Li, H. J. Zhu :}
 Zeta functions of totally ramified $p$-covers of the projective line.
{\it Rend. Sem. Mat. Univ. Padova,}
{\bf 113} (2005), 203--225.

\bibitem{Liu}
{\sc C. Liu :}
Generic exponential sums associated to Laurent polynomials in one variable, Preprint, 2008.

\bibitem{Ro}
{\sc P. Robba :} Index of $p$-adic differential operators III. Applications to twisted exponential sums.
{\it Ast\'erisque,}
{\bf 119-120}  (1984), 191--266.

\bibitem{Wan1} 
{\sc D. Wan :} Newton polygons for zeta and $L$-functions, {\it Ann. Math.} {\bf 137} (1993), 249--296.

\bibitem{Wan2} 
{\sc D. Wan :}
Variation of $p$-adic Newton polygons for $L$-functions of exponential sums, {\it Asian J. Math.} {\bf 8} (2004), 427--472.

\bibitem{Zhu2}
{\sc H. J. Zhu:}
$p$-adic variation of $L$ functions of one variable exponential
sums, I. {\it Amer. J. Math.} {\bf 125} (2003).

\bibitem{Zhu1}
{\sc H. J. Zhu :}
Asymptotic variation of $L$-functions of one-variable exponential
sums.
{\it J. Reine Angew. Math.} {\bf 572} (2004),
219--233. 

\end{thebibliography}
\end{document}